# A socio-physics based hybrid metaheuristic for solving complex non-convex constrained optimization problems

**Ishaan R Kale[1], Anand J Kulkarni[1], Efren Mezura-Montes[2]**

[1]Institute of Artificial Intelligence, MIT World Peace University,
Pune 411 038, MH, INDIA
ishaan.kale@mitwpu.edu.in; kale.ishaan@gmail.com
anand.j.kulkarni@mitwpu.edu.in; kulk0003@ntu.edu.sg

[2]Artificial Intelligence Research Institute, University of Veracruz, MEXICO
emezura@uv.mx

**Abstract**

Several Artificial Intelligence based heuristic and metaheuristic algorithms have been developed so far. These algorithms have shown their superiority towards solving complex problems from different domains. However, it is necessary to critically validate these algorithms for solving real-world constrained optimization problems. The search behavior in those problems is different as it involves large number of linear, nonlinear and non-convex type equality and inequality constraints. In this work a 57 real-world constrained optimization problems test suite is solved using two constrained metaheuristic algorithms originated from a socio-based Cohort Intelligence (CI) algorithm. The first CI-based algorithm incorporates a self-adaptive penalty function approach i.e., CI-SAPF. The second algorithm combines CI-SAPF with the intrinsic properties of the physics-based Colliding Bodies Optimization (CBO) referred to CI-SAPF-CBO. The results obtained from CI-SAPF and CI-SAPF-CBO are compared with other constrained optimization algorithms such as IUDE, $\epsilon$MAg-ES and iLSHADE$\epsilon$. The superiority of the proposed algorithms is discussed in details followed by future directions to evolve the constrained handling techniques.

## 1. Introduction

Many nature-inspired algorithms have been developed within the framework of Artificial Intelligence (AI) such as Genetic Algorithm (GA) (Deb and Goyal, 1996; Wu and Chow, 1995), Particle Swarm Optimization (PSO) (Li et al., 2009), Ant Colony Optimization (ACO) (Dorigo et al., 2006), Firefly Algorithm (FA) (Gandomi et al., 2011), Harmony Search Algorithm (HCA) (Lee et al., 2005), etc. to solve the complex search problems. These nature-inspired algorithms can be sub-categorized into bio-inspired, swarm, socio/cultural-based, inspired from laws of physics and chemistry (Kale and Kulkarni, 2021). These algorithms were initially designed to solve problems from a specific domain and are referred as *heuristic* algorithms. Considering this as one of the limitations, and to deal with problems from different domains such as engineering, health care and logistics, economics, etc., *metaheuristic* algorithms (Talbi, 2009) were proposed. From the literature, it has been noticed that these algorithms / problem solving methodologies requires certain modifications or need to hybridized with other suitable technique in order to improve its problems solving ability. It should be noted that the performance of modified or hybrid algorithm must be better or comparable as compared to standalone algorithm.

Generally, real-world optimization problems are associated with continuous, discrete and mixed variables as well as linear, nonlinear and nonconvex type equality and inequality constraints. Such elements generate different sources of difficulty regarding the constraints that the search algorithm has to del with. In order to consider feasibility information, an appropriate constraint-handling technique is required. There are several constraint handling techniques discussed in the literature such as the penalty function approach, feasibility-based rules (Deb, 2000, Kulkarni et al., 2016), probability-based constraint-handling technique (Kulkarni and Shabir, 2016a), epsilon constraint-handling (Takahama et al., 2005), dominance-based selection scheme (Coello and Mezura-Montes, 2002), decoders, special operators and separation of objective function and constraints (Mezura-Montes and Coello, 2011). The most popular constraint-handling technique is the penalty function (Homaifar et al., 1994) which converts the constrained problem into an unconstrained problem. However, it is driven by penalty parameters which need to be arbitrarily chosen to penalize infeasible solutions. Such issue gets complicated when several constraints are involved. This can lead to poor quality of the solution reached. To overcome such limitations, several penalty function variants have been developed such as dynamic penalty function (Joines and Houck, 1994), death penalty function (Luenberger and Ye, 2016), annealing penalty function (Michalewicz and Attia, 1994; Carlson et al., 1998), niched-penalty approach (Deb and Agrawal, 1999), self-adaptive penalty approach (Coello, 2000). However, these penalty-based constraint-handling techniques require external parameter/s which may lead to fine-tuning. This process could be computationally expensive as it needs several preliminary trials.

In the specialized literature, several nature-inspired algorithms are incorporated with different constraint-handling approaches discussed above to solve complex constrained optimization problems. Similar to other nature-inspired algorithms, Cohort Intelligence (CI) is one of those socio-inspired optimization algorithms (Kulkarni et al., 2013). CI is inspired from the self-supervised learning candidates having inherently common goal to achieve the best possible behavior which helps to evolve the behavior of the entire cohort. Several constrained and unconstrained optimization problems have been solved using the CI algorithm. It is important to note that CI was also incorporated with different constraint-handling technique such as the probability-based constraint-handling approach, static, dynamic and self-adaptive penalty function (SAPF) approaches. So far, CI with these penalty function approaches has been applied to solve problems from truss structural domain, engineering design (Kale and Kulkarni, 2018, 2021) and manufacturing domains (Kale et al., 2022). The CI with static penalty function (CI-SPF) approach required to set a suitable penalty parameter to penalize infeasible solutions. Hence, the SAPF approach was proposed, which does not require parameter tuning. CI was also applied to solve the 0-1 Knapsack problem (Kulkarni and Shabir, 2016a), the travelling salesman problem (Kulkarni et al., 2017), a healthcare and inventory management problem, a sea-cargo mix problem and a cross-border shipper selection (Kulkarni et al., 2016b). The constraints involved with these problems were highly nonlinear and interdependent. A probability-based constraint-handling technique was used to deal with the constraints. This constraint-handling technique is problem-specific and tedious to apply when the number of constraints increases. On the other hand, in Kale and Kulkarni (2021), the CI algorithm was hybridized with the physics-based Colliding Bodies Optimization (CBO) (i.e., CI-SAPF-CBO) algorithm to avoid the sampling parameter setting. Further, the CI-SAPF-CBO mechanism improved the exploration and exploitation of the search space.

In the present work, real-world constrained optimization problems are chosen to assess the performance of the CI-SAPF and CI-SAPF-CBO algorithms (Kale and Kulkarni, 2021). These problems are from six different domains such as (i) industrial and chemical process, (ii) process design and synthesis, (iii) mechanical design engineering, (iv) power system, (v) power electronics and (vi) livestock feed ration optimization (Kumar et al. 2020). This constrained optimization test suit was previously solved using three constrained evolutionary algorithms: Unified Differential Evolution (IUDE) algorithm (Trivedi et al., 2018), Modified Matrix Adaptation Evolution Strategy with Restarts ($\epsilon$MAg-ES) (Hellwig and Beyer, 2018) and an improved constrained version of an algorithm having IEpsilon, Linear Population Size Reduction and Success-history based parameter adaptation with Differential Evolution (iLSHADE$\epsilon$) (Fan et al., 2018). The results from CI-SAPF and CI-SAPF-CBO are compared with IUDE, $\epsilon$MAg-ES and iLSHADE$\epsilon$.

The manuscript is organized as follows: The CI algorithm is explained in Section 2. The framework of CI-SAPF along with its pseudo code are presented in Section 2.1. In Section 2.2, the framework of a hybridized CI-SAPF-CBO algorithm along with its pseudo code are presented. Section 3 gives an overview about the real-world constrained optimization test suit and presents the result comparison of CI-SAPF and CI-SAPF-CBO with other contemporary algorithms. The result discussion and analysis are presented in Section 4. Section 5 concludes the work and end up with the future directions.

## 2. Cohort Intelligence (CI)

Cohort Intelligence is the self-supervised socio-based Artificial Intelligence algorithm (Kulkarni et al., 2013). It models the behavior of learning candidates $C$ having inherently common goal to achieve the best possible behavior. Every candidate $c$ $(c = 1,2, \dots, C)$ in the cohort learns, follows and competes with every other candidate to evolve its individual behavior. The candidates in the cohort employs a probability-based roulette wheel approach to follow a candidate in the cohort which results to improve its own behavior. This makes every candidate learn from one another and helps the entire cohort behavior to evolve. The cohort behavior could be considered saturated, if for considerable number of learning attempts the behavior of every candidate does not considerably improve and become almost the same.

### 2.1 Framework of CI-SAPF

In general, the constrained optimization problem is expressed as follows:

Minimize $f(\boldsymbol{X}) = f(X_1, X_2, \dots, X_i, \dots, X_n)$ (2.1)

Subject to

$$g_i(\boldsymbol{X}) \leq 0, \qquad i = 1, 2, \dots, p$$
$$h_i(\boldsymbol{X}) = 0, \qquad i = 1, 2, \dots, m$$
$$\boldsymbol{\Psi}^{lower} \leq \boldsymbol{X} \leq \boldsymbol{\Psi}^{upper}$$

In a cohort with $C$ number of candidates. For every individual candidate $c(c = 1,2, \dots, C)$, using the CI-SAPF approach the pseudo-objective function (behavior) (refer to Eq. 2.1) can be expressed as follows:

$$\phi(\boldsymbol{X^c}) = f(\boldsymbol{X^c}) + SAPF(\boldsymbol{X^c}) \qquad (2.2)$$

where $SAPF(\boldsymbol{X^c}) = f(\boldsymbol{X^c}) \times (\sum_{i=1}^{n} g_i(\boldsymbol{X^c}) + \sum_{i=1}^{m} h_i(\boldsymbol{X^c}))$ is penalty function and $f(\boldsymbol{X^c})$ is the objective function of individual candidate. The behavior of the penalty parameter is dependent on the variable sampling space because, as the algorithm progresses it narrows down the sampling space using a sampling space reduction factor $R$. An independent penalty parameter $f(\boldsymbol{X^c})$ is generated by every individual candidate $c$ to penalize the violated constraints associated with its behavior. And further updates the penalty parameter subsequently for every learning attempt of CI-SAPF algorithm. However, the limitation of CI-SAPF has been observed. At the end of every learning attempt (iteration) the cohort candidates update its individual search space using a sampling space reduction factor $R$. The choice of $R$ is decided based on preliminary trials. To overcome this limitation, CI-SAPF is hybridized with CBO. It is discussed in section 2.2.

The above Eq. (2.1) is only applicable for positive function value i.e., $f(\boldsymbol{X})$. There are various other applications (numerical problems) having negative function value i.e., $-f(\boldsymbol{X})$ for those $SAPF$ is expressed as follows:

$$SAPF = abs(f(\boldsymbol{X})) \times \left(\sum_{i=1}^{p} g_i(\boldsymbol{X}) + \sum_{i=1}^{m} h_i(\boldsymbol{X})\right) \qquad (2.3)$$

And further forms the pseudo objective function $\phi(\boldsymbol{X})$ as follows:

$$\phi(\boldsymbol{X}) = abs(-f(\boldsymbol{X}) + SAPF). \qquad (2.4)$$

Specifically, when function value of any of the problem is too small or zero the $SAPF$ approach may not give feasible solution. In such cases, an arbitrary integer (fixed) value $(int)$ is added in the function value $f(\boldsymbol{X})$. Then the $SAPF$ would be calculated as follows:

$$SAPF = (f(\boldsymbol{X}) + int) \times (\sum_{i=1}^{n} g_i(\boldsymbol{X}) + \sum_{i=1}^{m} h_i(\boldsymbol{X})) \qquad (2.5)$$

Whereas, for the function obtaining infinity value $f(\boldsymbol{X}) = \infty$ the computation would be stop. In such cases, in order to continue the computation $f(\boldsymbol{X})$ need to be considered as unity or non-zero number. The pseudo code for CI-SAPF is presented in Figure 1.

| | |
|---|---|
| $C$ | Number of candidates $(c = 1, \dots, C)$ |
| $\boldsymbol{\Psi}$ | Sampling space |
| $R$ | Sampling space reduction factor |
| $P$ | Penalty function |
| $g$ | Violated constraint value |

Initialize $C, \boldsymbol{\Psi}, R$

**While**

1. Generate the arbitrary solutions using uniform distribution method.
2. Evaluate individual behavior/objective function $f(\mathbf{X^c})$
3. Apply SAPF approach to generate a pseudo behavior/objective function:
   $\phi(\mathbf{X^c}) = f(\mathbf{X^c}) + P(\mathbf{X^c})$
   where $P(\mathbf{X^c}) = f(\mathbf{X^c}) \times \sum_{c=1}^{C} g\,(\mathbf{X^c})$
4. The probability $p^c$ associated with every candidate $c$ in the cohort is calculated as: $p^c = \frac{1/\phi^*(\mathbf{X}^c)}{\sum_{c=1}^{C} 1/\phi^*(\mathbf{X}^c)}$
5. Using roulette wheel approach every candidate $c$ selects behavior to follow from within $C$ available choices.
6. Every candidate $c$ shrinks the sampling interval $\boldsymbol{\Psi^c}$ in its neighborhood using sampling space reduction parameter $R$:
   $\left[\boldsymbol{\Psi}^{c,lower}, \boldsymbol{\Psi}^{c,upper}\right] = \left[\mathbf{X}^c - \left\|\frac{\boldsymbol{\Psi}^{upper} - \boldsymbol{\Psi}^{lower}}{2}\right\| \times R,\ \mathbf{X}^c + \left\|\frac{\boldsymbol{\Psi}^{upper} - \boldsymbol{\Psi}^{lower}}{2}\right\| \times R\right]$
7. **If**: No significant improvement in the behavior $\phi^*(\mathbf{X}^c)$ is observed the solution is considered to be saturated
   Every candidate $c$ expands the sampling interval $\boldsymbol{\Psi^c}$ to its original interval $\boldsymbol{\Psi}$
   Accept the current behavior of cohort $\phi(\mathbf{X})$ and the associated attributes $\mathbf{X}$.
   **Else**
   Continue to Step 3
   **End If**

Figure 1: Pseudo Code of CI-SAPF algorithm (Kale and Kulkarni, 2021)

### 2.2 Hybrid CI-SAPF-CBO

The CI-SAPF-CBO is a hybrid metaheuristic algorithm models the intrinsic characteristics of CI and CBO algorithms (Kale and Kulkarni, 2021). As mentioned earlier, the performance of CI algorithm is dependent on a sampling space reduction factor ($R$) which require to be set for every problem. The motivation behind to propose a CI-SAPF-CBO is to improve the exploration rate of search space and to avoid the tuning of sampling space reduction factor $R$ (refer Step 6 in Figure 1). The CI-SAPF-CBO algorithm employs CI for global search, SAPF for constraint handling and CBO for the local search. The natural tendency of CI candidates is to follow candidates chosen probabilistically using roulette wheel approach to evolve individual behavior. For the sake of understanding, the terminology of CBO algorithm is expresses in the context of CI-SAPF algorithm i.e., stationary bodies, moving bodies, velocity and collision are expresses as slow learning candidates, fast learning candidates, learning ability and follow, respectively. Further, the learning ability of CI candidates is refined (updated) using CBO. It means that the fast-learning candidates motivates the slow learning candidate to improve the position in the search space. Further, the Coefficient of Restitution (COR) control the exploration of the search space and exploitation of the best solution. It also controls the local and global search. The pseudo code of CI-SAPF-CBO is presented in Figure 2.

| | |
|---|---|
| $C$ | Number of candidates ($c = 1, \dots, C$) |
| $\boldsymbol{\Psi}$ | Sampling space |
| $R$ | Sampling space reduction factor |
| $P$ | Penalty function |
| $g$ | Violated constraint value |

Initialize $C, \boldsymbol{\Psi}, R$

**While**

1. Generate the arbitrary solutions using uniform distribution method.
2. Evaluate individual behavior/objective function $f(\mathbf{X}^c)$
3. Apply SAPF approach to generate a pseudo behavior/objective function:
   $\phi(\mathbf{X}^c) = f(\mathbf{X}^c) + P(\mathbf{X}^c)$
   where $P(\mathbf{X}^c) = f(\mathbf{X}^c) \times \sum_{c=1}^{C} g(\mathbf{X}^c)$
4. The probability $p^c$ associated with every candidate $c$ in the cohort is calculated as: $p^c = \frac{1/\phi^*(\mathbf{X}^c)}{\sum_{c=1}^{C} 1/\phi^*(\mathbf{X}^c)}$
5. Using roulette wheel approach every candidate $c$ selects behavior to follow from within $C$ available choices.
6. First half of the candidates are referred as stationary bodies (slow learning candidates) and second half of candidates are referred as moving bodies (fast learning candidates).
7. Set the velocity (learning ability) of the slow learning candidates before follow
   $v^c = 0, \quad c = 1,2, \dots, C/2$
8. Evaluate learning ability of fast learning candidates before follow
   $v^c = \boldsymbol{X}^c - \boldsymbol{X}^{c-\frac{C}{2}}, \quad c = C/2 + 1, C/2 + 2, \dots, C$
9. Evaluate learning ability of fast learning candidates after follow
   $v'^c = \frac{\left(m^c - \varepsilon m^{c-C/2}\right) v^c}{m^c + m^{c-C/2}}, \quad c = C/2 + 1, C/2 + 2, \dots, C$
10. Evaluate learning ability of slow learning candidates after follow
    $m'^c = \frac{\left(p^{c+C/2} + \varepsilon p^{c+C/2}\right) v^{c+C/2}}{p^c + p^{c+C/2}}, \quad c = 1,2, \dots, C/2$
11. Evaluate the position of slow learning candidates $\boldsymbol{X}^{c_{new}} = \boldsymbol{X}^c + rand.\, v'^c, \quad c = 1,2, \dots, C/2$
12. Evaluate the position of fast learning candidates
    $\boldsymbol{X}^{c_{new}} = \boldsymbol{X}^{c-C/2} + rand.\, v'^c, \quad c = C/2 + 1, C/2 + 2, \dots, C$
13. Continue to step 2

**End**

Figure 2: Pseudo Code of the CI-SAPF-CBO algorithm (Kale and Kulkarni, 2021)

## 3 Solution to Real World Constrained Optimization Test Suit

In this paper, 57 real-world constrained optimization problems (Kumar et al., 2020) are considered to investigate the applicability of CI-SAPF and CI-SAPF-CBO (Kale and Kulkarni, 2021). These benchmarks constrained test problems are having problems from different domains such as industrial and chemical process, process design and synthesis, mechanical design engineering, power system, power electronics and livestock feed ration optimization. They are consisting of nonlinear function as well as nonlinear nonconvex equality and inequality constraints. The generalized optimization problem formulation is as follows:

Minimize $\quad F(\boldsymbol{X}), \boldsymbol{X} = (x_1, x_2, \dots, x_v)$

Subject to $\quad g_i(\boldsymbol{X}) \leq 0, i = 1, \dots, m$

$$h_j(\boldsymbol{X}) = 0, j = 1, \dots, n$$

where, $F$ is the objective function, $g$ is the inequality constraint and $h$ is the equality constraint. Equality constraint $h$ can be transformed into two inequality constraints using the following equation,

$$|h_j(\boldsymbol{X})| - \varepsilon \leq 0, j = 1, \ldots, n$$

where, $\varepsilon$ is set to a small value $10^{-4}$.

The simulations were run on Windows platform using an Intel(R) Core (TM)2Duo, 2.93GHz processor speed and 4GB RAM. Furthermore, every individual problem was solved 30 times. The solutions obtained from proposed techniques and comparison with other contemporary algorithms are discussed in the following sections.

Table 1: Comparison of results with the best function value available in the literature
($D$-number of design variables, $g$-number of inequality constraint and $h$-number of equality constraints)

| Problem | Name | $D$ | $g$ | $h$ | Best function value available in literature | CI-SAPF | CI-SAPF-CBO |
|---|---|---|---|---|---|---|---|
| **Industrial Chemical Processes** | | | | | | | |
| 1 | Heat Exchanger Network Design (case 1) | 9 | 0 | 8 | 1.8931162966E+02 | Infeasible | Infeasible |
| 2 | Heat Exchanger Network Design (case 2) | 11 | 0 | 9 | 7.0490369540E+03 | Infeasible | Infeasible |
| 3 | Optimal Operation of Alkylation Unit | 7 | 14 | 0 | -4.5291197395E+03 | -1999.9563 | -1995.5072 |
| 4 | Reactor Network Design (RND) | 6 | 1 | 4 | −3.8826043623E-01 | -0.9914 | -0.3711 |
| 5 | Haverly's Pooling Problem | 9 | 2 | 4 | −4.0000560000E+02 | -370.3079 | -225.0000 |
| 6 | Blending-Pooling-Separation problem | 38 | 0 | 32 | 1.8638304088E+00 | 1.0041 | 1.4118 |
| 7 | Propane, Isobutane, n-Butane Nonsharp Separation | 48 | 0 | 38 | 2.1158627569E+00 | 1.1179 | 1.1775 |
| **Process Synthesis and Design Problems** | | | | | | | |
| 8 | Process synthesis problem | 2 | 2 | 0 | 2.0000000000E+00 | 2.0000 | 2.0000 |
| 9 | Process synthesis and design problem | 3 | 1 | 1 | 2.5576545740E+00 | 1.5000 | 1.5000 |
| 10 | Process flow sheeting problem | 3 | 3 | 0 | 1.0765430833E+00 | 1.0765 | 1.0780 |
| 11 | Two-reactor Problem | 7 | 4 | 4 | 9.9238463653E+1 | Infeasible | Infeasible |
| 12 | Process synthesis problem | 7 | 9 | 0 | 2.9248305537E+00 | 2.9291 | 2.9591 |
| 13 | Process design Problem | 5 | 3 | 0 | 2.6887000000E+04 | 26887.4270 | 26899.0192 |
| 14 | Multi-product batch plant | 10 | 10 | 0 | 5.3638942722E+04 | 57213.4297 | 56817.9333 |
| **Mechanical Engineering Problems** | | | | | | | |
| 15 | Weight Minimization of a Speed Reducer | 7 | 11 | 0 | 2.9944244658E+03 | 2817.0905 | 2816.3410 |
| 16 | Optimal Design of Industrial refrigeration System | 14 | 15 | 0 | 3.2213000814E-02 | Infeasible | Infeasible |
| 17 | Tension/compression spring design (case 1) | 3 | 3 | 0 | 1.2665232788E-02 | 0.01272 | 0.01269 |
| 18 | Pressure vessel design | 4 | 4 | 0 | 5.8853327736E+03 | 5850.6640 | 5850.6640 |
| 19 | Welded beam design | 4 | 5 | 0 | 1.6702177263E+00 | 1.5511 | 1.5492 |
| 20 | Three-bar truss design problem | 2 | 3 | 0 | 2.6389584338E+02 | 263.8959 | 263.8959 |
| 21 | Multiple disk clutch brake design problem | 5 | 6 | 0 | 2.3524245790E-01 | 0.2352 | 0.2352 |
| 22 | Planetary gear train design optimization problem | 9 | 10 | 1 | 5.2576870748E-01 | 0.5258 | 0.5258 |
| 23 | Step-cone pulley problem | 5 | 8 | 3 | 1.6069868725E+01 | 14.2828 | 15.8163 |
| 24 | Robot gripper problem | 7 | 7 | 0 | 2.5287918415E+00 | 3.5002 | 4.1384 |
| 25 | Hydro-static thrust bearing design problem | 4 | 7 | 0 | 1.6254428092E+03 | Infeasible | Infeasible |
| 26 | Four-stage gear box problem | 22 | 86 | 0 | 3.5359231973E+01 | 41.6039 | 36.5620 |
| 27 | 10-bar truss design | 10 | 3 | 0 | 5.2445076066E+02 | 524.5624 | 528.9794 |
| 28 | Rolling element bearing | 10 | 9 | 0 | 1.4614135715E+04 | 17099.0558 | 14688.4578 |
| 29 | Gas Transmission Compressor Design (GTCD) | 4 | 1 | 0 | 2.9648954173E+06 | 2966066.0454 | 3168144.6401 |
| 30 | Tension/compression spring design (case 2) | 3 | 8 | 0 | 2.6138840583E+00 | 2.6586 | 2.6586 |
| 31 | Gear train design Problem | 4 | 1 | 1 | 0.0000000000E+00 | 2.7009E-12 | 2.7009E-12 |
| 32 | Himmelblau's Function | 5 | 6 | 0 | −3.0665538672E+04 | -30665.5384 | -30512.4500 |
| 33 | Topology Optimization | 30 | 30 | 0 | 2.6393464970E+00 | 2.7610 | 2.7899 |
| **Power System Problems** | | | | | | | |
| 34 | Optimal Sizing of Single Phase Distributed Generation with reactive power support for Phase Balancing at Main Transformer/Grid | 118 | 0 | 108 | 0.0000000000E+00 | Infeasible | Infeasible |
| 35 | Optimal Sizing of Distributed Generation for Active Power Loss Minimization | 153 | 0 | 148 | 8.9093896456E-02 | Infeasible | Infeasible |
| 36 | Optimal Sizing of Distributed Generation (DG) and Capacitors for Reactive Power Loss Minimization | 158 | 0 | 148 | 7.2066551720E-02 | Infeasible | Infeasible |
| 37 | Optimal Power flow (Minimization of Active Power Loss) | 126 | 0 | 116 | 2.1962851478E-02 | Infeasible | Infeasible |
| 38 | Optimal Power flow (Minimization of Fuel Cost) | 126 | 0 | 116 | 2.7766131989E+00 | Infeasible | Infeasible |
| 39 | Optimal Power flow (Minimization of Active Power Loss and Fuel Cost) | 126 | 0 | 116 | 2.8677165770E+00 | Infeasible | Infeasible |
| 40 | Microgrid Power flow (Islanded case) | 76 | 0 | 76 | 0.0000000000E+00 | Infeasible | Infeasible |
| 41 | Microgrid Power flow (Grid-connected case) | 74 | 0 | 74 | 0.0000000000E+00 | Infeasible | Infeasible |
| 42 | Optimal Setting of Droop Controller for Minimization of Active Power Loss in Islanded Microgrids | 86 | 0 | 86 | 8.6241006360E-02 | Infeasible | Infeasible |
| 43 | Optimal Setting of Droop Controller for Minimization of Reactive Power Loss in Islanded Microgrids | 86 | 0 | 76 | 8.0420545897E-02 | Infeasible | Infeasible |
| 44 | Wind Farm Layout Problem | 3 | 91 | 0 | −6.2607000000E+03 | -5664.3528 | -5591.2511 |
| **Power Electronic Problems** | | | | | | | |
| 45 | SOPWM for 3-level Inverters | 25 | 24 | 1 | 3.8029250566E-02 | 0.0004 | Infeasible |
| 46 | SOPWM for 5-level Inverters | 25 | 24 | 1 | 2.1215000000E-02 | 0.0371 | Infeasible |
| 47 | SOPWM for 7-level Inverters | 25 | 24 | 1 | 1.5164538375E-02 | 0.0199 | Infeasible |
| 48 | SOPWM for 9-level Inverters | 30 | 29 | 1 | 1.6787535766E-02 | 0.0395 | Infeasible |
| 49 | SOPWM for 11-level Inverters | 30 | 29 | 1 | 9.3118741800E-03 | 0.0172 | Infeasible |
| 50 | SOPWM for 13-level Inverters | 30 | 29 | 1 | 1.5096451396E-02 | 0.0141 | Infeasible |
| **Livestock Feed Ration Optimization** | | | | | | | |
| 51 | Beef Cattle (case 1) | 59 | 14 | 1 | 4.5508511497E+03 | Infeasible | Infeasible |
| 52 | Beef Cattle (case 2) | 59 | 14 | 1 | 3.3489821493E+03 | Infeasible | Infeasible |
| 53 | Beef Cattle (case 3) | 59 | 14 | 1 | 4.9976069290E+03 | Infeasible | Infeasible |
| 54 | Beef Cattle (case 4) | 59 | 14 | 1 | 4.2405482538E+03 | Infeasible | Infeasible |
| 55 | Dairy Cattle (case 1) | 64 | 0 | 6 | 6.6964145128E+03 | Infeasible | Infeasible |
| 56 | Dairy Cattle (case 2) | 64 | 0 | 6 | 1.4748932529E+04 | Infeasible | Infeasible |
| 57 | Dairy Cattle (case 3) | 64 | 0 | 6 | 3.2132917019E+03 | Infeasible | Infeasible |

Table 2: Comparison of statistical results obtained using CI-SAPF and CI-SAPF-CBO

| Problem | Algorithm | Best | Median | Mean | Worst | STD | MCV | FR | Avg. FE | Avg. CT $(sec)$ |
|---|---|---|---|---|---|---|---|---|---|---|
| **Industrial Chemical Process Problems** | | | | | | | | | | |
| 1 | IUDE | 1.89E+02 | 2.60E+02 | 2.29E+02 | 1.85E+02 | 8.06E+01 | 1.12E+05 | 24 | NA | NA |
| | $\epsilon$MAg-ES | 1.89E+02 | 4.92E+02 | 4.55E+02 | 4.37E+02 | 4.37E+02 | 1.47E-04 | 84 | NA | NA |
| | iLSHADE$\epsilon$ | 1.90E+02 | 1.94E+02 | 2.06E+02 | 2.29E+02 | 1.93E+01 | 1.36E-02 | 28 | NA | NA |
| | **CI-SAPF** | Infeasible | Infeasible | Infeasible | Infeasible | Infeasible | Best 2.24E+06<br>Mean 6.33E+06<br>Worst 1.15E+07 | - | - | - |
| | **CI-SAPF-CBO** | Infeasible | Infeasible | Infeasible | Infeasible | Infeasible | Best 2.27E+06<br>Mean 7.10E+06<br>Worst 1.24E+07 | - | - | - |
| 2 | IUDE | 7.05E+03 | 7.05E+03 | 7.15E+03 | 5.94E+03 | 7.54E+02 | 6.17E+03 | 92 | NA | NA |
| | $\epsilon$MAg-ES | 7.05E+03 | 7.80E+03 | 7.74E+03 | 7.48E+03 | 7.50E+02 | 3.71E+03 | 96 | NA | NA |
| | iLSHADE$\epsilon$ | 7.05E+03 | 7.05E+03 | 7.05E+03 | 7.05E+03 | 2.13E-11 | 0.00E+00 | 100 | NA | NA |
| | **CI-SAPF** | Infeasible | Infeasible | Infeasible | Infeasible | Infeasible | Best 7.10E+06<br>Mean 7.42E+06<br>Worst 7.73E+06 | - | - | - |
| | **CI-SAPF-CBO** | Infeasible | Infeasible | Infeasible | Infeasible | Infeasible | Best 3.83E+06<br>Mean 6.47E+06<br>Worst 7.61E+06 | - | - | - |
| 3 | IUDE | −4.53E+03 | −1.43E+02 | −6.25E+03 | −1.83E+04 | 6.75E+03 | 3.66E+00 | 64 | NA | NA |
| | $\epsilon$MAg-ES | −1.43E+02 | 7.63E+01 | −1.60E+02 | 4.95E+02 | 8.88E+02 | 1.61E-02 | 92 | NA | NA |
| | iLSHADE$\epsilon$ | −4.53E+03 | −1.43E+02 | −8.18E+02 | 5.34E+02 | 1.91E+03 | 0.00E+00 | 100 | NA | NA |
| | **CI-SAPF** | **-1999.9563** | **-1999.9924** | **-1999.9887** | **-1999.9999** | **0.0122** | **0.00E+00** | **100** | **4976** | **137.99** |
| | **CI-SAPF-CBO** | **-1995.5073** | **-1999.8882** | **-1999.5923** | **-1999.9669** | **1.1318** | **0.00E+00** | **100** | **5162** | **125.52** |
| 4 | IUDE | −2.86E-01 | −5.92E-01 | −4.96E-01 | −1.00E+00 | 1.82E-0 | 7.90E-02 | 0 | NA | NA |
| | $\epsilon$MAg-ES | −3.88E-01 | −3.75E-01 | −5.50E-01 | −1.00E+00 | 2.86E-01 | 3.50E-02 | 72 | NA | NA |
| | iLSHADE$\epsilon$ | −3.75E-01 | −3.75E-01 | −3.75E-01 | −3.73E-01 | 4.61E-04 | 0.00E+00 | 100 | NA | NA |
| | **CI-SAPF** | **-0.9999** | **-0.9998** | **-0.9989** | **-0.9914** | **0.0021** | **0.00E+00** | **100** | **2571** | **4.32** |
| | **CI-SAPF-CBO** | **-1** | **-0.4106** | **-0.49142** | **-0.3711** | **0.2075** | **0.00E+00** | **100** | **348** | **2.84** |
| 5 | IUDE | −4.00E+02 | −4.00E+02 | −3.51E+02 | −8.30E-03 | 1.32E+02 | 0.00E+00 | 100 | NA | NA |
| | $\epsilon$MAg-ES | −4.00E+02 | −3.98E+02 | −3.63E+02 | −1.00E+02 | 7.55E+01 | 0.00E+00 | 100 | NA | NA |
| | iLSHADE$\epsilon$ | −4.00E+02 | −8.06E-03 | −1.17E+02 | 1.57E+01 | 1.79E+02 | 0.00E+00 | 100 | NA | NA |
| | **CI-SAPF** | **-370.3079** | **-586.4450** | **-573.3752** | **-836.1784** | **125.8807** | **0.00E+00** | **100** | **4038** | **111.98** |
| | **CI-SAPF-CBO** | **-225.0000** | **-2.77E+02** | **-306.4025** | **-529.0002** | **96.0494** | **0.00E+00** | **100** | **2216** | **18.77** |
| 6 | IUDE | 1.71E+00 | 9.98E-01 | 1.07E+00 | 9.98E-01 | 1.98E-01 | 2.10E+00 | 0 | NA | NA |
| | $\epsilon$MAg-ES | 2.09E+00 | 1.83E+00 | 1.59E+00 | 1.20E+00 | 3.39E-01 | 7.07E-02 | 0 | NA | NA |
| | iLSHADE$\epsilon$ | 1.72E+00 | 1.24E+00 | 1.27E+00 | 1.11E+00 | 1.63E-01 | 3.94E+00 | 0 | NA | NA |
| | **CI-SAPF** | **1.0041** | **1.0682** | **1.0964** | **1.2658** | **0.0811** | **0.00E+00** | **100** | **3201** | **88.77** |
| | **CI-SAPF-CBO** | **1.4118** | **2.6957** | **2.6463** | **4.2711** | **0.9631** | **0.00E+00** | **100** | **4835** | **157** |
| 7 | IUDE | 1.76E+00 | 1.34E+00 | 1.40E+00 | 9.98E-01 | 3.82E-01 | 2.82E-01 | 0 | NA | NA |
| | $\epsilon$MAg-ES | 2.00E+00 | 1.71E+00 | 1.80E+00 | 2.01E+00 | 1.74E-01 | 1.98E-02 | 0 | NA | NA |
| | iLSHADE$\epsilon$ | 1.75E+00 | 1.73E+00 | 1.66E+00 | 1.67E+00 | 1.00E-01 | 2.99E+00 | 0 | NA | NA |
| | **CI-SAPF** | **1.1179** | **1.2875** | **1.3684** | **1.8154** | **0.2513** | **0.00E+00** | **100** | **4057** | **136.45** |
| | **CI-SAPF-CBO** | **1.1775** | **1.5861** | **1.8912** | **3.237498718** | **0.667773597** | **0.00E+00** | **100** | **5432** | **44.28** |
| **Process Synthesis and Design Problems** | | | | | | | | | | |
| 8 | IUDE | 2.00E+00 | 2.00E+00 | 2.00E+00 | 2.00E+00 | 6.41E-17 | 0.00E+00 | 100 | NA | NA |
| | $\epsilon$MAg-ES | 2.00E+00 | 2.00E+00 | 1.99E+00 | 1.29E+00 | 1.52E-01 | 4.58E-03 | 96 | NA | NA |
| | iLSHADE$\epsilon$ | 2.00E+00 | 2.00E+00 | 2.00E+00 | 2.00E+00 | 0.00E+00 | 0.00E+00 | 100 | NA | NA |
| | **CI-SAPF** | **2.0000** | **2.0000** | **2.0945** | **2.2361** | **0.1197** | **0.00E+00** | **100** | **1244** | **5.99** |
| | **CI-SAPF-CBO** | **2.0000** | **2.0000** | **2.0787** | **2.2361** | **0.1152** | **0.00E+00** | **100** | | **5.52** |
| 9 | IUDE | 2.56E+00 | 2.56E+00 | 2.56E+00 | 2.56E+00 | 1.36E-15 | 0.00E+00 | 100 | NA | NA |
| | $\epsilon$MAg-ES | 2.56E+00 | 2.56E+00 | 2.55E+00 | 1.93E+00 | 2.70E-01 | 1.15E-02 | 92 | NA | NA |
| | iLSHADE$\epsilon$ | 2.56E+00 | 2.56E+00 | 2.56E+00 | 2.56E+00 | 1.46E-07 | 0.00E+00 | 100 | NA | NA |
| | **CI-SAPF** | **1.5000** | **1.5000** | **1.5000** | **1.5000** | **1.2867E-05** | **0.00E+00** | **100** | **1089** | **4.20** |
| | **CI-SAPF-CBO** | **1.5000** | **1.5000** | **1.5000** | **1.5000** | **1.2272E-05** | **0.00E+00** | **100** | **3810** | **5.44** |
| 10 | IUDE | 1.08E+00 | 1.08E+00 | 1.11E+00 | 1.25E+00 | 7.08E-02 | 0.00E+00 | 100 | NA | NA |
| | $\epsilon$MAg-ES | 1.08E+00 | 1.08E+00 | 1.08E+00 | 1.25E+00 | 3.47E-02 | 0.00E+00 | 100 | NA | NA |
| | iLSHADE$\epsilon$ | 1.08E+00 | 1.25E+00 | 1.22E+00 | 1.25E+00 | 6.48E-02 | 0.00E+00 | 100 | NA | NA |
| | **CI-SAPF** | **1.0765** | **1.0765** | **1.0765** | **1.0765** | **4.3957E-07** | **0.00E+00** | **100** | **3741** | **30.57** |
| | **CI-SAPF-CBO** | **1.0780** | **1.0854** | **1.1186** | **1.2501** | **0.0677** | **0.00E+00** | **100** | **9364** | **22.56** |
| 11 | IUDE | 9.92E+01 | 9.92E+01 | 1.02E+02 | 1.07E+02 | 4.07E+00 | 0.00E+00 | 100 | NA | NA |
| | $\epsilon$MAg-ES | 1.07E+02 | 9.92E+01 | 1.06E+02 | 1.15E+02 | 6.88E+00 | 2.20E-02 | 0 | NA | NA |
| | iLSHADE$\epsilon$ | 9.92E+01 | 1.07E+02 | 1.06E+02 | 1.32E+02 | 7.39E+00 | 4.00E-02 | 96 | NA | NA |
| | **CI-SAPF** | Infeasible | Infeasible | Infeasible | Infeasible | Infeasible | Best 4.82E+01<br>Mean 5.88E+01<br>Worst 7.22E+01 | - | - | - |
| | **CI-SAPF-CBO** | Infeasible | Infeasible | Infeasible | Infeasible | Infeasible | Best 4.31E+01<br>Mean 5.14E+01<br>Worst 7.71E+01 | - | - | - |
| 12 | IUDE | 2.92E+00 | 2.95E+00 | 3.00E+00 | 4.21E+00 | 2.54E-01 | 0.00E+00 | 100 | NA | NA |
| | $\epsilon$MAg-ES | 2.92E+00 | 3.64E+00 | 3.65E+00 | 4.69E+00 | 5.87E-01 | 0.00E+00 | 100 | NA | NA |
| | iLSHADE$\epsilon$ | 2.92E+00 | 2.92E+00 | 2.92E+00 | 2.92E+00 | 8.14E-07 | 0.00E+00 | 100 | NA | NA |
| | **CI-SAPF** | **2.9291** | **2.9559** | **2.9667** | **3.0692** | **0.0382** | **0.00E+00** | **100** | **2859** | **62.54** |
| | **CI-SAPF-CBO** | **2.9591** | **3.0077** | **3.0736** | **3.4156** | **0.1357** | **0.00E+00** | **100** | **13431** | **22.42** |
| 13 | IUDE | 2.69E+04 | 2.69E+04 | 2.69E+04 | 2.69E+04 | 1.11E-11 | 0.00E+00 | 100 | NA | NA |
| | $\epsilon$MAg-ES | 2.69E+04 | 2.69E+04 | 2.69E+04 | 2.69E+04 | 1.11E-11 | 0.00E+00 | 100 | NA | NA |
| | iLSHADE$\epsilon$ | 2.69E+04 | 2.69E+04 | 2.69E+04 | 2.69E+04 | 1.11E-11 | 0.00E+00 | 100 | NA | NA |
| | **CI-SAPF** | **26887.4270** | **26887.9139** | **26897.7311** | **26945.0523** | **20.1523** | **0.00E+00** | **100** | **2560** | **20.47** |
| | **CI-SAPF-CBO** | **26899.0192** | **26921.3667** | **26950.4005** | **27277.5399** | **92.9462** | **0.00E+00** | **100** | **15276** | **47.50** |
| 14 | IUDE | 5.85E+04 | 6.65E+04 | 6.60E+04 | 7.36E+04 | 5.14E+03 | 0.00E+00 | 100 | NA | NA |
| | $\epsilon$MAg-ES | 5.36E+04 | 5.85E+04 | 5.81E+04 | 5.85E+04 | 1.35E+03 | 0.00E+00 | 100 | NA | NA |
| | iLSHADE$\epsilon$ | 5.36E+04 | 5.92E+04 | 5.91E+04 | 6.36E+04 | 1.85E+03 | 0.00E+00 | 100 | NA | NA |
| | **CI-SAPF** | **57213.4297** | **62088.7485** | **61968.6733** | **67501.4561** | **3408.6197** | **0.00E+00** | **100** | **1920** | **26.06** |
| | **CI-SAPF-CBO** | **56817.9333** | **63694.2212** | **63070.0756** | **66568.4494** | **2686.5258** | **0.00E+00** | **100** | **24489** | **32.89** |
| **Mechanical Engineering Problems** | | | | | | | | | | |
| 15 | IUDE | 2.99E+03 | 2.99E+03 | 2.99E+03 | 2.99E+03 | 4.64E-13 | 0.00E+00 | 100 | NA | NA |
| | $\epsilon$MAg-ES | 2.99E+03 | 2.99E+03 | 2.99E+03 | 2.99E+03 | 4.64E-13 | 0.00E+00 | 100 | NA | NA |
| | iLSHADE$\epsilon$ | 2.99E+03 | 2.99E+03 | 2.99E+03 | 2.99E+03 | 4.64E-13 | 0.00E+00 | 100 | NA | NA |
| | **CI-SAPF** | **2817.0905** | **2820.5046** | **2820.3191** | **2825.8713** | **2.0660** | **0.00E+00** | **100** | **2060** | **2.25** |
| | **CI-SAPF-CBO** | **2816.3410** | **2817.4050** | **2817.6529** | **2819.5744** | **1.0276** | **0.00E+00** | **100** | **1946** | **3.61** |
| 16 | IUDE | 3.22E-02 | 3.22E-02 | 6.24E+00 | 2.95E+00 | 2.86E+01 | 1.28E-02 | 80 | NA | NA |
| | $\epsilon$MAg-ES | 3.22E-02 | 3.22E-02 | 3.22E-02 | 3.22E-02 | 2.78E-17 | 0.00E+00 | 100 | NA | NA |
| | iLSHADE$\epsilon$ | 3.22E-02 | 3.23E-02 | 3.82E-01 | 2.95E+00 | 9.67E-01 | 1.16E-01 | 80 | NA | NA |
| | **CI-SAPF** | Infeasible | Infeasible | Infeasible | Infeasible | Infeasible | Best 6.97E+00<br>Mean 1.05E+01<br>Worst 1.47E+01 | - | - | - |
| | **CI-SAPF-CBO** | Infeasible | Infeasible | Infeasible | Infeasible | Infeasible | Best 2.29E-01<br>Mean 1.30E+01<br>Worst 4.39E+01 | - | - | - |
| 17 | IUDE | 1.27E-02 | 1.27E-02 | 1.27E-02 | 1.27E-02 | 1.08E-05 | 0.00E+00 | 100 | NA | NA |
| | $\epsilon$MAg-ES | 1.27E-02 | 1.27E-02 | 1.27E-02 | 1.37E-02 | 2.16E-04 | 0.00E+00 | 100 | NA | NA |
| | iLSHADE$\epsilon$ | 1.27E-02 | 1.27E-02 | 1.30E-02 | 1.78E-02 | 1.06E-03 | 0.00E+00 | 100 | NA | NA |
| | **CI-SAPF** | **0.0127** | **0.01311** | **0.01380** | **0.01729** | **0.0015** | **0.00E+00** | **100** | **1295** | **8.52** |
| | **CI-SAPF-CBO** | **0.0126** | **0.01297** | **0.01308** | **0.01473** | **0.0004** | **0.00E+00** | **100** | **2964** | **20.96** |
| 18 | IUDE | 6.06E+03 | 6.06E+03 | 6.06E+03 | 6.09E+03 | 6.16E+00 | 0.00E+00 | 100 | NA | NA |
| | $\epsilon$MAg-ES | 6.06E+03 | 6.41E+03 | 7.38E+03 | 1.19E+04 | 1.93E+03 | 0.00E+00 | 100 | NA | NA |
| | iLSHADE$\epsilon$ | 6.06E+03 | 6.11E+03 | 8.48E+03 | 1.49E+04 | 3.14E+03 | 0.00E+00 | 100 | NA | NA |
| | **CI-SAPF** | **5850.6640** | **5916.9771** | **5960.4924** | **6116.8912** | **104.6492** | **0.00E+00** | **100** | **2155** | **4.85** |
| | **CI-SAPF-CBO** | **5850.6640** | **5895.6063** | **5946.4602** | **6095.0539** | **97.0569** | **0.00E+00** | **100** | **2443** | **4.69** |
| 19 | IUDE | 1.67E+00 | 1.67E+00 | 1.67E+00 | 1.67E+00 | 1.20E-16 | 0.00E+00 | 100 | NA | NA |
| | $\epsilon$MAg-ES | 1.67E+00 | 1.67E+00 | 1.69E+00 | 1.85E+00 | 3.95E-02 | 0.00E+00 | 100 | NA | NA |
| | iLSHADE$\epsilon$ | 1.67E+00 | 1.67E+00 | 1.67E+00 | 1.67E+00 | 7.59E-07 | 0.00E+00 | 100 | NA | NA |
| | **CI-SAPF** | **1.5511** | **1.5603** | **1.5604** | **1.5742** | **0.0077** | **0.00E+00** | **100** | **2874** | **4.79** |
| | **CI-SAPF-CBO** | **1.5492** | **1.5502** | **1.5527** | **1.5591** | **0.0035** | **0.00E+00** | **100** | **2180** | **2.24** |
| 20 | IUDE | 2.64E+02 | 2.64E+02 | 2.64E+02 | 2.64E+02 | 0.0000 | 0.00E+00 | 100 | NA | NA |
| | ϵMAg-ES | 2.64E+02 | 2.64E+02 | 2.65E+02 | 2.74E+02 | 2.88E+00 | 0.00E+00 | 100 | NA | NA |
| | iLSHADEϵ | 2.64E+02 | 2.64E+02 | 2.64E+02 | 2.64E+02 | 1.99E-02 | 0.00E+00 | 100 | NA | NA |
| | **CI-SAPF** | **263.8959** | **263.8963** | **263.8972** | **263.9064** | **0.0024** | **0.00E+00** | **100** | **2378** | **16.14** |
| | **CI-SAPF-CBO** | **263.8959** | **263.8981** | **263.8994** | **263.9128** | **0.0050** | **0.00E+00** | **100** | **2740** | **16.10** |
| 21 | IUDE | 2.35E-01 | 2.35E-01 | 2.35E-01 | 2.35E-01 | 1.13E-16 | 0.00E+00 | 100 | NA | NA |
| | ϵMAg-ES | 2.35E-01 | 2.35E-01 | 2.35E-01 | 2.35E-01 | 1.13E-16 | 0.00E+00 | 100 | NA | NA |
| | iLSHADEϵ | 2.35E-01 | 2.35E-01 | 2.35E-01 | 2.35E-01 | 1.13E-16 | 0.00E+00 | 100 | NA | NA |
| | **CI-SAPF** | **0.235242458** | **0.235242458** | **0.235242458** | **0.235242458** | **8.61893E-17** | **0.00E+00** | **100** | **165** | **0.45** |
| | **CI-SAPF-CBO** | **0.235242458** | **0.235242458** | **0.235242458** | **0.235242458** | **8.61893E-17** | **0.00E+00** | **100** | **190** | **0.44** |
| 22 | IUDE | 5.26E-01 | 5.26E-01 | 5.26E-01 | 5.27E-01 | 5.50E-04 | 0.00E+00 | 100 | NA | NA |
| | ϵMAg-ES | 5.26E-01 | 5.46E-01 | 6.16E-01 | 1.12E+00 | 1.98E-01 | 1.31E+00 | 76 | NA | NA |
| | iLSHADEϵ | 5.26E-01 | 5.26E-01 | 5.27E-01 | 5.31E-01 | 1.69E-03 | 0.00E+00 | 100 | NA | NA |
| | **CI-SAPF** | **0.5258** | **0.5287** | **0.5291** | **0.5421** | **0.0037** | **0.00E+00** | **100** | **1104** | **10.87** |
| | **CI-SAPF-CBO** | **0.5258** | **0.5306** | **0.5336** | **0.5258** | **0.0072** | **0.00E+00** | **100** | **1114** | **15.00** |
| 23 | IUDE | 1.61E+01 | 1.61E+01 | 1.61E+01 | 1.61E+01 | 4.17E-15 | 0.00E+00 | 100 | NA | NA |
| | ϵMAg-ES | 1.61E+01 | 1.61E+01 | 1.61E+01 | 1.61E+01 | 1.78E-14 | 0.00E+00 | 100 | NA | NA |
| | iLSHADEϵ | 1.61E+01 | 1.61E+01 | 1.61E+01 | 1.61E+01 | 8.62E-08 | 0.00E+00 | 100 | NA | NA |
| | **CI-SAPF** | **14.2828** | **14.2983** | **14.3084** | **14.3804** | **0.0302** | **0.00E+00** | **100** | **3794** | **1.72** |
| | **CI-SAPF-CBO** | **15.8163** | **16.1601** | **16.1974** | **17.3040** | **0.4316** | **0.00E+00** | **100** | **2691** | **2.56** |
| 24 | IUDE | 2.54E+00 | 2.54E+00 | 2.54E+00 | 2.54E+00 | 5.99E-14 | 0.00E+00 | 100 | NA | NA |
| | ϵMAg-ES | 2.54E+00 | 2.54E+00 | 2.54E+00 | 2.55E+00 | 8.81E-04 | 0.00E+00 | 100 | NA | NA |
| | iLSHADEϵ | 2.54E+00 | 2.54E+00 | 2.54E+00 | 2.55E+00 | 1.99E-03 | 0.00E+00 | 100 | NA | NA |
| | **CI-SAPF** | **3.5002** | **4.6432** | **4.8001** | **6.3996** | **0.7119** | **0.00E+00** | **100** | **1500** | **5.59** |
| | **CI-SAPF-CBO** | **4.1384** | **4.8134** | **5.0750** | **7.1009** | **0.8012** | **0.00E+00** | **100** | **1793** | **5.57** |
| 25 | IUDE | 1.86E+03 | 2.60E+02 | 1.93E+03 | 2.60E+02 | 2.37E+03 | 9.01E-05 | 40 | NA | NA |
| | ϵMAg-ES | 1.62E+03 | 2.26E+03 | 2.35E+03 | 6.34E+02 | 1.41E+03 | 1.84E-05 | 88 | NA | NA |
| | iLSHADEϵ | 1.66E+03 | 2.09E+03 | 1.76E+03 | 6.12E+02 | 1.28E+03 | 2.54E-04 | 76 | NA | NA |
| | **CI-SAPF** | Convergence not achieved | | | | | | | | |
| | **CI-SAPF-CBO** | Convergence not achieved | | | | | | | | |
| 26 | IUDE | 4.54E+01 | 4.91E+01 | 6.21E+01 | 4.55E+01 | 2.93E+01 | 1.35E-01 | 12 | NA | NA |
| | ϵMAg-ES | 8.23E+01 | 2.13E+01 | 3.38E+01 | 7.26E+00 | 3.61E+01 | 3.77E-01 | 8 | NA | NA |
| | iLSHADEϵ | 3.54E+01 | 3.63E+01 | 3.63E+01 | 3.73E+01 | 6.59E-01 | 0.00E+00 | 100 | NA | NA |
| | **CI-SAPF** | **41.6039** | **50.8290** | **55.7953** | **96.3675** | **13.8099** | **0.00E+00** | **100** | **2532** | **15.53** |
| | **CI-SAPF-CBO** | **36.5620** | **64.7603** | **80.1487** | **128.6801** | **29.2113** | **0.00E+00** | **100** | **833** | **14.69** |
| 27 | IUDE | 5.24E+02 | 5.24E+02 | 5.24E+02 | 5.24E+02 | 4.72E-04 | 0.00E+00 | 100 | NA | NA |
| | ϵMAg-ES | 5.24E+02 | 5.31E+02 | 5.30E+02 | 5.31E+02 | 2.03E+00 | 0.00E+00 | 100 | NA | NA |
| | iLSHADEϵ | 5.24E+02 | 5.25E+02 | 5.25E+02 | 5.25E+02 | 7.14E-02 | 0.00E+00 | 100 | NA | NA |
| | **CI-SAPF** | **524.5624** | **524.9151** | **525.2325** | **527.3440** | **0.7366** | **0.00E+00** | **100** | **3040** | **12.37** |

| | | | | | | | | | | | |
|---|---|---|---|---|---|---|---|---|---|---|---|
| | CI-SAPF-CBO | 528.9794 | 533.5394 | 533.9862 | 539.5492 | 3.6486 | 0.00E+00 | | 100 | 3648 | 15.63 |
| 28 | IUDE | 1.46E+04 | 1.46E+04 | 1.46E+04 | 1.46E+04 | 9.28E-12 | 0.00E+00 | | 100 | NA | NA |
| | ϵMAg-ES | 1.46E+04 | 1.46E+04 | 1.46E+04 | 1.46E+04 | 9.28E-12 | 0.00E+00 | | 100 | NA | NA |
| | iLSHADEϵ | 1.46E+04 | 1.46E+04 | 1.46E+04 | 1.46E+04 | 9.28E-12 | 0.00E+00 | | 100 | NA | NA |
| | CI-SAPF | 17099.0558 | 17280.8071 | 17309.7824 | 17726.2415 | 188.1210 | 0.00E+00 | | 100 | 1959 | 24.24 |
| | CI-SAPF-CBO | 14688.4578 | 15334.7574 | 15518.9090 | 17671.4797 | 838.0081 | 0.00E+00 | | 100 | 1943 | 17.97 |
| 29 | IUDE | 2.96E+06 | 2.96E+06 | 2.96E+06 | 2.96E+06 | 1.43E-09 | 0.00E+00 | | 100 | NA | NA |
| | ϵMAg-ES | 2.96E+06 | 2.96E+06 | 2.96E+06 | 2.96E+06 | 1.43E-09 | 0.00E+00 | | 100 | NA | NA |
| | iLSHADEϵ | 2.96E+06 | 2.96E+06 | 2.97E+06 | 2.97E+06 | 6.57E+02 | 0.00E+00 | | 100 | NA | NA |
| | CI-SAPF | 2966066.0454 | 2969337.5205 | 2984136.3990 | 3072153.3581 | 25237.5419 | 0.00E+00 | | 100 | 1527 | 11.01 |
| | CI-SAPF-CBO | 3168144.6401 | 3168144.6401 | 3168144.6401 | 3168144.6401 | 0.0000 | 0.00E+00 | | 100 | 663 | 4.53 |
| 30 | IUDE | 2.66E+00 | 4.54E+00 | 5.53E+00 | 2.67E+01 | 4.74E+00 | 1.81E-03 | | 92 | NA | NA |
| | ϵMAg-ES | 2.66E+00 | 3.11E+00 | 2.21E+00 | 4.34E-02 | 1.22E+00 | 1.04E+06 | | 68 | NA | NA |
| | iLSHADEϵ | 6.90E+00 | 2.87E+00 | 6.26E+00 | 2.05E+01 | 6.32E+00 | 7.95E-01 | | 20 | NA | NA |
| | CI-SAPF | 2.6586 | 2.6586 | 2.6613 | 2.7034 | 0.0096 | 0.00E+00 | | 100 | 3494 | 19.67 |
| | CI-SAPF-CBO | 2.6586 | 2.6588 | 2.6601 | 2.6722 | 0.0032 | 0.00E+00 | | 100 | 795 | 4.74 |
| 31 | IUDE | 0.00E+00 | 1.74E-18 | 4.55E-16 | 8.41E-15 | 1.68E-15 | 0.00E+00 | | 100 | NA | NA |
| | ϵMAg-ES | 0.00E+00 | 0.00E+00 | 0.00E+00 | 0.00E+00 | 0.00E+00 | 0.00E+00 | | 100 | NA | NA |
| | iLSHADEϵ | 0.00E+00 | 5.41E-18 | 5.56E-17 | 3.91E-16 | 1.17E-16 | 0.00E+00 | | 100 | NA | NA |
| | CI-SAPF | 2.7009E-12 | 2.7009E-12 | 4.05939E-12 | 2.30782E-11 | 5.2614E-12 | 0.00E+00 | | 100 | 2815 | 3.12 |
| | CI-SAPF-CBO | 2.7009E-12 | 2.31E-11 | 1.72561E-11 | 2.30782E-11 | 9.55302E-12 | 0.00E+00 | | 100 | 1498 | 2.05 |
| 32 | IUDE | −3.07E+04 | −3.07E+04 | −3.07E+04 | −3.07E+04 | 3.71E-12 | 0.00E+00 | | 100 | NA | NA |
| | ϵMAg-ES | −3.07E+04 | −3.07E+04 | −3.07E+04 | −3.07E+04 | 3.56E-12 | 0.00E+00 | | 100 | NA | NA |
| | iLSHADEϵ | −3.07E+04 | −3.07E+04 | −3.07E+04 | −3.07E+04 | 3.64E-12 | 0.00E+00 | | 100 | NA | NA |
| | CI-SAPF | -30665.5384 | -30665.1049 | -30659.9610 | -30630.6691 | 9.6236 | 0.00E+00 | | 100 | 5957 | 50.44 |
| | CI-SAPF-CBO | -30512.4500 | -30452.9387 | -30451.6920 | -30512.4500 | 44.1197 | 0.00E+00 | | 100 | 639 | 6.38 |
| 33 | IUDE | 2.64E+00 | 2.64E+00 | 2.64E+00 | 2.64E+00 | 1.41E-15 | 0.00E+00 | | 100 | NA | NA |
| | ϵMAg-ES | 2.64E+00 | 2.65E+00 | 2.65E+00 | 2.68E+00 | 1.26E-02 | 0.00E+00 | | 40 | NA | NA |
| | iLSHADEϵ | 2.64E+00 | 2.64E+00 | 2.64E+00 | 2.64E+00 | 8.11E-16 | 0.00E+00 | | 100 | NA | NA |
| | CI-SAPF | 2.7610 | 2.8570 | 2.8671 | 3.0918 | 0.0845 | 0.00E+00 | | 100 | 1472 | 14.71 |
| | CI-SAPF-CBO | 2.7899 | 2.9509 | 2.9686 | 3.2947 | 0.1071 | 0.00E+00 | | 100 | 1477 | 14.77 |
| **Power System Problems** | | | | | | | | | | | |
| 34 | IUDE | 5.07E+00 | 1.17E+01 | 5.01E+00 | 2.38E+00 | 1.99E+00 | 4.27E-02 | | 0 | NA | NA |
| | ϵMAg-ES | 3.06E+00 | 3.75E+00 | 5.50E+00 | 6.53E+00 | 2.87E+00 | 1.43E-02 | | 0 | NA | NA |
| | iLSHADEϵ | 1.04E+01 | 6.73E+00 | 8.85E+00 | 5.76E+00 | 2.49E+00 | 5.10E+00 | | 0 | NA | NA |
| | CI-SAPF | Infeasible | Infeasible | Infeasible | Infeasible | Infeasible | Best | 3.54E-01 | - | - | - |
| | | | | | | | Mean | 9.11E+00 | | | |
| | | | | | | | Worst | 2.67E+01 | | | |
| | CI-SAPF-CBO | Infeasible | Infeasible | Infeasible | Infeasible | Infeasible | Best | 7.56E-02 | - | - | - |
| | | | | | | | Mean | 1.42E+01 | | | |
| | | | | | | | Worst | 6.05E+01 | | | |
| 35 | IUDE | 9.33E+01 | 9.68E+01 | 1.01E+02 | 1.03E+02 | 1.32E+01 | 8.52E-01 | | 0 | NA | NA |
| | $\epsilon$MAg-ES | 8.80E+01 | 4.68E+01 | 7.77E+01 | 9.40E+01 | 2.09E+01 | 9.24E-01 | | 0 | NA | NA |
| | iLSHADE$\epsilon$ | 1.88E+02 | 1.87E+02 | 1.57E+02 | 1.64E+02 | 2.49E+01 | 2.75E+01 | | 0 | NA | NA |
| | **CI-SAPF** | Infeasible | Infeasible | Infeasible | Infeasible | Infeasible | Best | 2.01E+02 | - | - | - |
| | | | | | | | Mean | 8.12E+02 | | | |
| | | | | | | | Worst | 1.36E+03 | | | |
| | **CI-SAPF-CBO** | Infeasible | Infeasible | Infeasible | Infeasible | Infeasible | Best | 5.73E+01 | - | - | - |
| | | | | | | | Mean | 8.27E+02 | | | |
| | | | | | | | Worst | 3.08E+03 | | | |
| 36 | IUDE | 8.08E+01 | 6.38E+01 | 8.16E+01 | 9.99E+01 | 1.82E+01 | 8.17E-01 | | 0 | NA | NA |
| | $\epsilon$MAg-ES | 7.69E+01 | 7.37E+01 | 7.70E+01 | 8.42E+01 | 8.96E+00 | 9.33E-01 | | 0 | NA | NA |
| | iLSHADE$\epsilon$ | 1.18E+02 | 9.62E+01 | 1.29E+02 | 1.36E+02 | 1.84E+01 | 5.77E+01 | | 0 | NA | NA |
| | **CI-SAPF** | Infeasible | Infeasible | Infeasible | Infeasible | Infeasible | Best | 2.17E+00 | - | - | - |
| | | | | | | | Mean | 9.69E+02 | | | |
| | | | | | | | Worst | 2.27E+03 | | | |
| | **CI-SAPF-CBO** | Infeasible | Infeasible | Infeasible | Infeasible | Infeasible | Best | 4.41E+00 | - | - | - |
| | | | | | | | Mean | 1.05E+03 | | | |
| | | | | | | | Worst | 2.27E+03 | | | |
| 37 | IUDE | 1.91E+00 | −1.88E+00 | 3.42E-01 | −6.31E+00 | 2.20E+00 | 1.23E-01 | | 0 | NA | NA |
| | $\epsilon$MAg-ES | 1.35E+00 | 1.57E+00 | 1.55E+00 | 5.37E-01 | 4.62E-01 | 1.94E-02 | | 0 | NA | NA |
| | iLSHADE$\epsilon$ | 3.02E+00 | 3.03E+00 | 3.48E+00 | 3.99E+00 | 4.53E-01 | 6.32E+00 | | 0 | NA | NA |
| | **CI-SAPF** | Infeasible | Infeasible | Infeasible | Infeasible | Infeasible | Best | 7.46E-01 | - | - | - |
| | | | | | | | Mean | 2.09E+01 | | | |
| | | | | | | | Worst | 5.06E+01 | | | |
| | **CI-SAPF-CBO** | Infeasible | Infeasible | Infeasible | Infeasible | Infeasible | Best | 1.11E+00 | - | - | - |
| | | | | | | | Mean | 1.88E+01 | | | |
| | | | | | | | Worst | 4.17E+01 | | | |
| 38 | IUDE | 3.67E+00 | −1.57E+01 | −1.41E+01 | −2.85E+01 | 1.02E+01 | 1.98E-01 | | 0 | NA | NA |
| | $\epsilon$MAg-ES | 6.09E+00 | 3.92E+00 | 4.17E+00 | 5.65E+00 | 1.69E+00 | 2.99E-02 | | 0 | NA | NA |
| | iLSHADE$\epsilon$ | 3.86E+00 | 4.17E+00 | 3.12E+00 | 1.44E+00 | 1.19E+00 | 6.67E+00 | | 0 | NA | NA |
| | **CI-SAPF** | Infeasible | Infeasible | Infeasible | Infeasible | Infeasible | Best | 1.67E+00 | - | - | - |
| | | | | | | | Mean | 1.09E+01 | | | |
| | | | | | | | Worst | 3.14E+01 | | | |
| | **CI-SAPF-CBO** | Infeasible | Infeasible | Infeasible | Infeasible | Infeasible | Best | 2.16E+00 | - | - | - |
| | | | | | | | Mean | 1.51E+01 | | | |
| | | | | | | | Worst | 3.98E+01 | | | |
| 39 | IUDE | −2.99E+00 | −1.58E+01 | −2.24E+01 | −4.00E+01 | 1.41E+01 | 2.09E-01 | | 0 | NA | NA |
| | $\epsilon$MAg-ES | 1.81E+00 | 3.56E+00 | 3.92E+00 | 2.36E+00 | 1.50E+00 | 3.01E-02 | | 0 | NA | NA |
| | iLSHADE$\epsilon$ | 5.01E+00 | 1.07E+00 | 3.18E+00 | 4.07E-01 | 1.65E+00 | 6.91E+00 | | 0 | NA | NA |
| | | | | | | | Best | 9.56E-01 | | | |
| | **CI-SAPF** | Infeasible | Infeasible | Infeasible | Infeasible | Infeasible | Mean | 1.94E+01 | - | - | - |
| | | | | | | | Worst | 4.97E+01 | | | |
| | | | | | | | Best | 5.18E+00 | | | |
| | **CI-SAPF-CBO** | Infeasible | Infeasible | Infeasible | Infeasible | Infeasible | Mean | 1.41E+01 | - | - | - |
| | | | | | | | Worst | 2.82E+01 | | | |
| 40 | IUDE | 8.61E+01 | 3.82E+01 | 1.16E+02 | 1.33E+02 | 6.66E+01 | 1.35E+00 | | 0 | NA | NA |
| | $\epsilon$MAg-ES | 3.05E-01 | 2.16E+01 | 3.02E+01 | 5.12E+01 | 2.38E+01 | 2.44E-01 | | 0 | NA | NA |
| | iLSHADE$\epsilon$ | 1.14E+02 | 1.11E+02 | 1.89E+02 | 6.00E+02 | 1.15E+02 | 1.50E+02 | | 0 | NA | NA |
| | | | | | | | Best | 8.98E+01 | | | |
| | **CI-SAPF** | Infeasible | Infeasible | Infeasible | Infeasible | Infeasible | Mean | 9.74E+02 | - | - | - |
| | | | | | | | Worst | 4.29E+03 | | | |
| | | | | | | | Best | 8.71E+01 | | | |
| | **CI-SAPF-CBO** | Infeasible | Infeasible | Infeasible | Infeasible | Infeasible | Mean | 6.24E+02 | - | - | - |
| | | | | | | | Worst | 1.61E+03 | | | |
| 41 | IUDE | 5.06E+01 | 9.19E+01 | 9.14E+01 | 6.18E+02 | 1.20E+02 | 1.14E+00 | | 0 | NA | NA |
| | $\epsilon$MAg-ES | 4.08E-08 | 3.02E-03 | 6.60E+00 | 3.15E+01 | 1.45E+01 | 6.71E-02 | | 0 | NA | NA |
| | iLSHADE$\epsilon$ | 4.17E+01 | 2.26E+02 | 1.24E+02 | 4.75E+02 | 1.06E+02 | 1.13E+02 | | 0 | NA | NA |
| | | | | | | | Best | 2.86E+01 | | | |
| | **CI-SAPF** | Infeasible | Infeasible | Infeasible | Infeasible | Infeasible | Mean | 6.07E+02 | - | - | - |
| | | | | | | | Worst | 1.37E+03 | | | |
| | | | | | | | Best | 4.14E+01 | | | |
| | **CI-SAPF-CBO** | Infeasible | Infeasible | Infeasible | Infeasible | Infeasible | Mean | 9.25E+02 | - | - | - |
| | | | | | | | Worst | 2.11E+03 | | | |
| 42 | IUDE | −1.70E+00 | 2.25E+02 | −3.64E+01 | −1.01E+03 | 2.48E+02 | 5.25E+00 | | 0 | NA | NA |
| | $\epsilon$MAg-ES | −1.33E+00 | 6.49E+01 | 9.24E+01 | 7.77E+01 | 3.89E+01 | 8.03E-01 | | 0 | NA | NA |
| | iLSHADE$\epsilon$ | −1.00E+00 | −1.01E+00 | −1.06E+00 | −1.17E+00 | 1.72E-01 | 1.65E+02 | | 0 | NA | NA |
| | | | | | | | Best | 3.55E+00 | | | |
| | **CI-SAPF** | Infeasible | Infeasible | Infeasible | Infeasible | Infeasible | Mean | 7.70E+02 | - | - | - |
| | | | | | | | Worst | 2.63E+03 | | | |
| | | | | | | | Best | 1.19E+02 | | | |
| | **CI-SAPF-CBO** | Infeasible | Infeasible | Infeasible | Infeasible | Infeasible | Mean | 8.61E+02 | - | - | - |
| | | | | | | | Worst | 2.84E+03 | | | |
| 43 | IUDE | 3.25E+01 | 1.18E+00 | 1.34E+01 | −1.73E+01 | 1.95E+01 | 2.53E+00 | | 0 | NA | NA |
| | $\epsilon$MAg-ES | 1.04E+02 | 1.18E+02 | 9.82E+01 | 9.85E+01 | 2.22E+01 | 8.21E-01 | | 0 | NA | NA |
| | iLSHADE$\epsilon$ | 2.72E+01 | 2.21E+01 | 3.80E+01 | 3.84E+01 | 7.86E+00 | 1.54E+02 | | 0 | NA | NA |
| | | | | | | | Best | 4.17E+01 | | | |
| | **CI-SAPF** | Infeasible | Infeasible | Infeasible | Infeasible | Infeasible | Mean | 6.38E+02 | - | - | - |
| | | | | | | | Worst | 1.70E+03 | | | |
| | **CI-SAPF-CBO** | Infeasible | Infeasible | Infeasible | Infeasible | Infeasible | Best | 8.58E+00 | **-** | - | - |
| | | | | | | | Mean | 7.43E+02 | | | |
| | | | | | | | Worst | 2.09E+03 | | | |
| 44 | IUDE | −6.09E+03 | −6.02E+03 | −6.00E+03 | −5.90E+03 | 4.61E+01 | 0.00E+00 | | 100 | NA | NA |
| | ϵMAg-ES | −6.17E+03 | −5.97E+03 | −5.97E+03 | −5.71E+03 | 1.05E+02 | 0.00E+00 | | 100 | NA | NA |
| | iLSHADEϵ | −6.34E+03 | −6.20E+03 | −6.20E+03 | −6.14E+03 | 3.66E+01 | 0.00E+00 | | 100 | NA | NA |
| | **CI-SAPF** | **-5664.3528** | **-5291.1577** | **-5240.2315** | **-4380.2667** | **348.4429** | **0.00E+00** | | **100** | **1899** | **18.99** |
| | **CI-SAPF-CBO** | **-5591.2511** | **-5300.8535** | **-5296.3133** | **-5030.0537** | **168.7838** | **0.00E+00** | | **100** | **1349** | **10.46** |
| **Power Electronics Problems** | | | | | | | | | | | |
| 45 | IUDE | 6.85E-02 | 1.13E-01 | 1.17E-01 | 2.49E-01 | 4.55E-02 | 0.00E+00 | | 100 | NA | NA |
| | ϵMAg-ES | 3.14E-02 | 4.94E-02 | 4.66E-02 | 5.71E-02 | 9.46E-03 | 0.00E+00 | | 100 | NA | NA |
| | iLSHADEϵ | 7.08E-02 | 1.51E-01 | 1.48E-01 | 2.49E-01 | 4.50E-02 | 0.00E+00 | | 100 | NA | NA |
| | **CI-SAPF** | **0.000455457** | **0.001581964** | **0.020827023** | **0.110157274** | **0.037032802** | **0.00E+00** | | **100** | **4397** | **15.01** |
| | **CI-SAPF-CBO** | Infeasible | Infeasible | Infeasible | Infeasible | Infeasible | Best | 3.95E-01 | **-** | - | - |
| | | | | | | | Mean | 2.61E+01 | | | |
| | | | | | | | Worst | 5.74E+01 | | | |
| 46 | IUDE | 2.02E-02 | 5.85E-02 | 5.54E-02 | 7.70E-02 | 1.58E-02 | 0.00E+00 | | 100 | NA | NA |
| | ϵMAg-ES | 2.04E-02 | 5.10E-02 | 4.82E-02 | 8.36E-02 | 1.25E-02 | 0.00E+00 | | 100 | NA | NA |
| | iLSHADEϵ | 5.17E-02 | 1.03E-01 | 1.05E-01 | 1.92E-01 | 3.75E-02 | 0.00E+00 | | 100 | NA | NA |
| | **CI-SAPF** | **0.03711278** | **0.09701265** | **0.147412781** | **0.464344792** | **0.111475993** | **0.00E+00** | | **100** | **5996** | **21.52** |
| | **CI-SAPF-CBO** | Infeasible | Infeasible | Infeasible | Infeasible | Infeasible | Best | 7.81E-01 | **-** | - | - |
| | | | | | | | Mean | 2.26E+01 | | | |
| | | | | | | | Worst | 8.11E+01 | | | |
| 47 | IUDE | 2.76E-02 | 5.37E-02 | 6.84E-02 | 2.06E-01 | 3.79E-02 | 0.00E+00 | | 100 | NA | NA |
| | ϵMAg-ES | 1.86E-02 | 3.19E-02 | 3.13E-02 | 5.16E-02 | 8.34E-03 | 0.00E+00 | | 100 | NA | NA |
| | iLSHADEϵ | 2.44E-02 | 5.10E-02 | 5.45E-02 | 1.01E-01 | 2.36E-02 | 0.00E+00 | | 100 | NA | NA |
| | **CI-SAPF** | **0.019908604** | **0.05153724** | **0.079316066** | **0.317820597** | **0.080113113** | **0.00E+00** | | **100** | **4663** | **16.56** |
| | **CI-SAPF-CBO** | Infeasible | Infeasible | Infeasible | Infeasible | Infeasible | Best | 2.36E+00 | **-** | - | - |
| | | | | | | | Mean | 3.79E+01 | | | |
| | | | | | | | Worst | 7.11E+01 | | | |
| 48 | IUDE | 4.64E-02 | 7.62E-02 | 1.00E-01 | 3.05E-01 | 6.43E-02 | 5.61E-04 | | 96 | NA | NA |
| | ϵMAg-ES | 2.24E-02 | 5.21E-02 | 6.14E-02 | 1.46E-01 | 2.89E-02 | 0.00E+00 | | 100 | NA | NA |
| | iLSHADEϵ | 4.64E-02 | 3.38E-01 | 2.92E-01 | 3.29E-01 | 1.50E-01 | 5.66E-02 | | 80 | NA | NA |
| | **CI-SAPF** | **0.039511817** | **0.142020919** | **0.167948863** | **0.55394813** | **0.127828637** | **0.00E+00** | | **100** | **4897** | **17.09** |
| | **CI-SAPF-CBO** | Infeasible | Infeasible | Infeasible | Infeasible | Infeasible | Best | 2.88E+00 | **-** | - | - |
| | | | | | | | Mean | 3.15E+01 | | | |
| | | | | | | | Worst | 8.23E+01 | | | |
| 49 | IUDE | 3.61E-02 | 6.96E-02 | 7.56E-02 | 2.84E-01 | 5.28E-02 | 0.00E+00 | | 100 | NA | NA |
| | ϵMAg-ES | 1.18E-02 | 2.47E-02 | 2.82E-02 | 5.51E-02 | 1.25E-02 | 0.00E+00 | | 100 | NA | NA |
| | iLSHADE$\epsilon$ | 6.74E-02 | 1.61E-01 | 1.85E-01 | 4.04E-01 | 7.26E-02 | 0.00E+00 | | 100 | NA | NA |
| | **CI-SAPF** | **0.017275557** | **0.045386079** | **0.045004006** | **0.104584644** | **0.021407327** | **0.00E+00** | | **100** | **7082** | **25.60** |
| | **CI-SAPF-CBO** | Infeasible | Infeasible | Infeasible | Infeasible | Infeasible | Best | 6.63E-01 | - | - | - |
| | | | | | | | Mean | 2.82E+01 | | | |
| | | | | | | | Worst | 7.12E+01 | | | |
| 50 | IUDE | 1.67E-01 | 3.34E-01 | 3.14E-01 | 1.70E-01 | 6.31E-02 | 6.05E-03 | | 16 | NA | NA |
| | $\epsilon$MAg-ES | 1.50E-02 | 1.62E-02 | 1.80E-02 | 3.51E-02 | 5.23E-03 | 0.00E+00 | | 100 | NA | NA |
| | iLSHADE$\epsilon$ | 2.66E-01 | 3.58E-01 | 3.32E-01 | 3.76E-01 | 4.49E-02 | 2.34E-01 | | 0 | NA | NA |
| | **CI-SAPF** | **0.014176355** | **0.033134891** | **0.037806182** | **0.058859149** | **0.01210712** | **0.00E+00** | | **100** | **6664** | **23.87** |
| | **CI-SAPF-CBO** | Infeasible | Infeasible | Infeasible | Infeasible | Infeasible | Best | 5.36E+00 | - | - | - |
| | | | | | | | Mean | 3.89E+01 | | | |
| | | | | | | | Worst | 7.57E+01 | | | |
| **Livestock Feed Ration Optimization** | | | | | | | | | | | |
| 51 | IUDE | 4.55E+03 | 4.55E+03 | 4.55E+03 | 4.55E+03 | 4.16E-02 | 2.82E-06 | | 0 | NA | NA |
| | $\epsilon$MAg-ES | 4.55E+03 | 4.55E+03 | 4.50E+03 | 3.05E+03 | 3.02E+02 | 1.67E-02 | | 0 | NA | NA |
| | iLSHADE$\epsilon$ | 4.55E+03 | 4.57E+03 | 4.56E+03 | 4.58E+03 | 1.14E+01 | 1.97E-04 | | 0 | NA | NA |
| | **CI-SAPF** | Infeasible | Infeasible | Infeasible | Infeasible | Infeasible | Best | 4.01E+00 | - | - | - |
| | | | | | | | Mean | 4.01E+00 | | | |
| | | | | | | | Worst | 4.01E+00 | | | |
| | **CI-SAPF-CBO** | Infeasible | Infeasible | Infeasible | Infeasible | Infeasible | Best | 4.01E+00 | - | - | - |
| | | | | | | | Mean | 4.01E+00 | | | |
| | | | | | | | Worst | 4.01E+00 | | | |
| 52 | IUDE | 3.35E+03 | 3.39E+03 | 3.39E+03 | 3.43E+03 | 2.16E+01 | 0.00E+00 | | 100 | NA | NA |
| | $\epsilon$MAg-ES | 3.45E+03 | 3.76E+03 | 3.82E+03 | 4.38E+03 | 2.47E+02 | 0.00E+00 | | 100 | NA | NA |
| | iLSHADE$\epsilon$ | 3.65E+03 | 3.96E+03 | 3.98E+03 | 4.28E+03 | 1.64E+02 | 0.00E+00 | | 100 | NA | NA |
| | **CI-SAPF** | Infeasible | Infeasible | Infeasible | Infeasible | Infeasible | Best | 1.39E+02 | - | - | - |
| | | | | | | | Mean | 1.77E+02 | | | |
| | | | | | | | Worst | 2.08E+02 | | | |
| | **CI-SAPF-CBO** | Infeasible | Infeasible | Infeasible | Infeasible | Infeasible | Best | 1.59E+02 | - | - | - |
| | | | | | | | Mean | 1.85E+02 | | | |
| | | | | | | | Worst | 2.18E+02 | | | |
| 53 | IUDE | 5.00E+03 | 5.07E+03 | 5.08E+03 | 5.31E+03 | 8.70E+01 | 0.00E+00 | | 100 | NA | NA |
| | $\epsilon$MAg-ES | 5.01E+03 | 5.25E+03 | 5.11E+03 | 2.57E+03 | 5.50E+02 | 3.56E-02 | | 92 | NA | NA |
| | iLSHADE$\epsilon$ | 5.06E+03 | 5.24E+03 | 5.27E+03 | 5.53E+03 | 1.56E+02 | 0.00E+00 | | 100 | NA | NA |
| | **CI-SAPF** | Infeasible | Infeasible | Infeasible | Infeasible | Infeasible | Best | 1.50E+02 | - | - | - |
| | | | | | | | Mean | 1.82E+02 | | | |
| | | | | | | | Worst | 2.15E+02 | | | |
| | **CI-SAPF-CBO** | Infeasible | Infeasible | Infeasible | Infeasible | Infeasible | Best | 1.51E+02 | - | - | - |
| | | | | | | | Mean | 1.89E+02 | | | |
| | | | | | | | Worst | 2.15E+02 | | | |
| 54 | IUDE | 4.24E+03 | 4.24E+03 | 4.24E+03 | 4.24E+03 | 2.10E-01 | 0.00E+00 | | 100 | NA | NA |
| | $\epsilon$MAg-ES | 4.24E+03 | 4.07E+03 | 3.69E+03 | 1.38E+03 | 7.31E+02 | 4.02E-02 | | 40 | NA | NA |
| | iLSHADE$\epsilon$ | 4.24E+03 | 4.24E+03 | 4.24E+03 | 4.25E+03 | 2.78E+00 | 0.00E+00 | | 100 | NA | NA |
| | **CI-SAPF** | Infeasible | Infeasible | Infeasible | Infeasible | Infeasible | Best | 1.63E+02 | - | - | - |
| | | | | | | | Mean | 2.00E+02 | | | |
| | | | | | | | Worst | 2.36E+02 | | | |
| | **CI-SAPF-CBO** | Infeasible | Infeasible | Infeasible | Infeasible | Infeasible | Best | 1.52E+02 | - | - | - |
| | | | | | | | Mean | 1.77E+02 | | | |
| | | | | | | | Worst | 2.10E+02 | | | |
| 55 | IUDE | 3.51E+03 | 2.06E+03 | 2.28E+03 | 2.23E+03 | 3.32E+02 | 9.71E-03 | | 0 | NA | NA |
| | $\epsilon$MAg-ES | 6.16E+03 | 2.67E+03 | 2.74E+03 | 3.61E+03 | 1.20E+03 | 6.84E-02 | | 0 | NA | NA |
| | iLSHADE$\epsilon$ | 6.27E+03 | 6.39E+03 | 5.38E+03 | 2.76E+03 | 1.20E+03 | 2.29E-02 | | 0 | NA | NA |
| | **CI-SAPF** | Infeasible | Infeasible | Infeasible | Infeasible | Infeasible | Best | 5.36E+01 | - | - | - |
| | | | | | | | Mean | 5.36E+01 | | | |
| | | | | | | | Worst | 5.36E+01 | | | |
| | **CI-SAPF-CBO** | Infeasible | Infeasible | Infeasible | Infeasible | Infeasible | Best | 5.36E+01 | - | - | - |
| | | | | | | | Mean | 5.36E+01 | | | |
| | | | | | | | Worst | 5.36E+01 | | | |
| 56 | IUDE | 1.41E+04 | 1.15E+04 | 1.19E+04 | 7.37E+03 | 1.49E+03 | 9.85E-03 | | 0 | NA | NA |
| | $\epsilon$MAg-ES | 1.32E+04 | 9.36E+03 | 9.12E+03 | 6.57E+03 | 2.23E+03 | 1.70E-02 | | 0 | NA | NA |
| | iLSHADE$\epsilon$ | 1.38E+04 | 1.44E+04 | 1.30E+04 | 1.12E+04 | 1.40E+03 | 5.16E-02 | | 0 | NA | NA |
| | **CI-SAPF** | Infeasible | Infeasible | Infeasible | Infeasible | Infeasible | Best | 4.75E+02 | - | - | - |
| | | | | | | | Mean | 6.08E+02 | | | |
| | | | | | | | Worst | 7.51E+02 | | | |
| | **CI-SAPF-CBO** | Infeasible | Infeasible | Infeasible | Infeasible | Infeasible | Best | 5.41E+02 | - | - | - |
| | | | | | | | Mean | 6.59E+02 | | | |
| | | | | | | | Worst | 8.28E+02 | | | |
| 57 | IUDE | 2.88E+03 | 2.72E+03 | 2.66E+03 | 2.39E+03 | 2.85E+02 | 1.98E-03 | | 0 | NA | NA |
| | $\epsilon$MAg-ES | 4.01E+03 | 2.61E+03 | 3.40E+03 | 1.06E+03 | 1.60E+03 | 2.34E-02 | | 0 | NA | NA |
| | iLSHADE$\epsilon$ | 8.31E+03 | 6.41E+03 | 6.98E+03 | 1.26E+04 | 2.62E+03 | 1.12E-02 | | 0 | NA | NA |
| | **CI-SAPF** | Infeasible | Infeasible | Infeasible | Infeasible | Infeasible | Best | 5.77E+02 | - | - | - |
| | | | | | | | Mean | 7.04E+02 | | | |
| | | | | | | | Worst | 8.22E+02 | | | |
| | **CI-SAPF-CBO** | Infeasible | Infeasible | Infeasible | Infeasible | Infeasible | Best | 5.50E+02 | - | - | - |
| | | | | | | | Mean | 6.58E+02 | | | |
| | | | | | | | Worst | 7.64E+02 | | | |

Table 3: Comparison of algorithms for finding feasible solution for number problems from different domains

| **Algorithms** | **Industrial Chemical Domain Problems** | **Process Synthesis and Design Domain problems** | **Mechanical engineering domain problems** | **Power system domain problems** | **Power electronic domain problems** | **Livestock feed ration optimization domain problems** | **Total feasibility out of 57** |
|---|---|---|---|---|---|---|---|
| IUDE | 1 | 7 | 15 | 1 | 4 | 2 | 30 |
| $\epsilon$MAg-ES | 1 | 4 | 16 | 1 | 6 | 1 | 29 |
| iLSHADE$\epsilon$ | 4 | 6 | 16 | 1 | 5 | 3 | 35 |
| **CI-SAPF** | 5 | 6 | 17 | 1 | 6 | 0 | 35 |
| **CI-SAPF-CBO** | 5 | 6 | 17 | 1 | 0 | 0 | 29 |

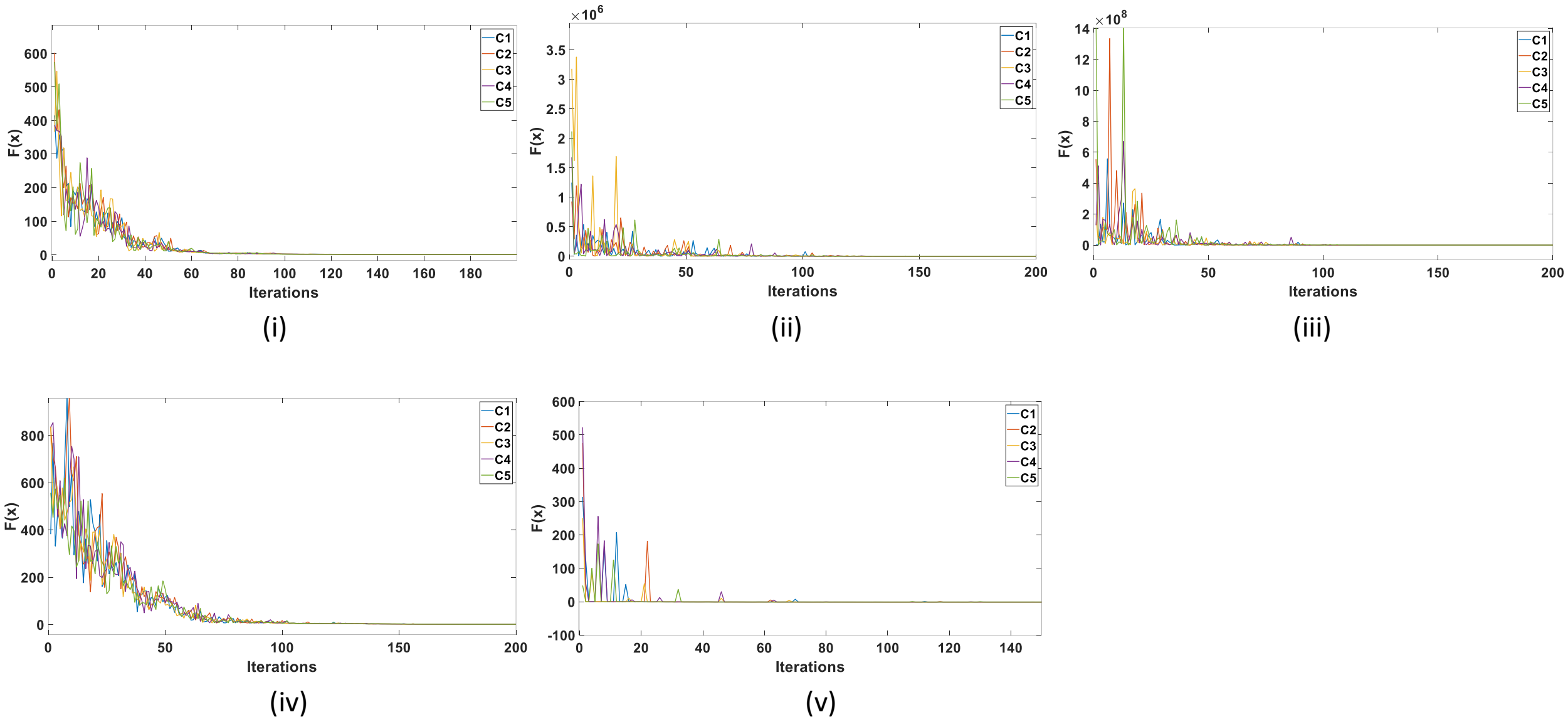

Fig. 3: Convergence plots of Industrial Chemical Process Problems using CI-SAPF algorithm
(i) Blending Pooling Separation Problem, (ii) Haverly's Pooling Problems, (iii) Optimal Operation of Alkylation Unit, (iv) Propane, Isobutane, n-Butane Nonsharp Separation, (v) Reactor Network Design

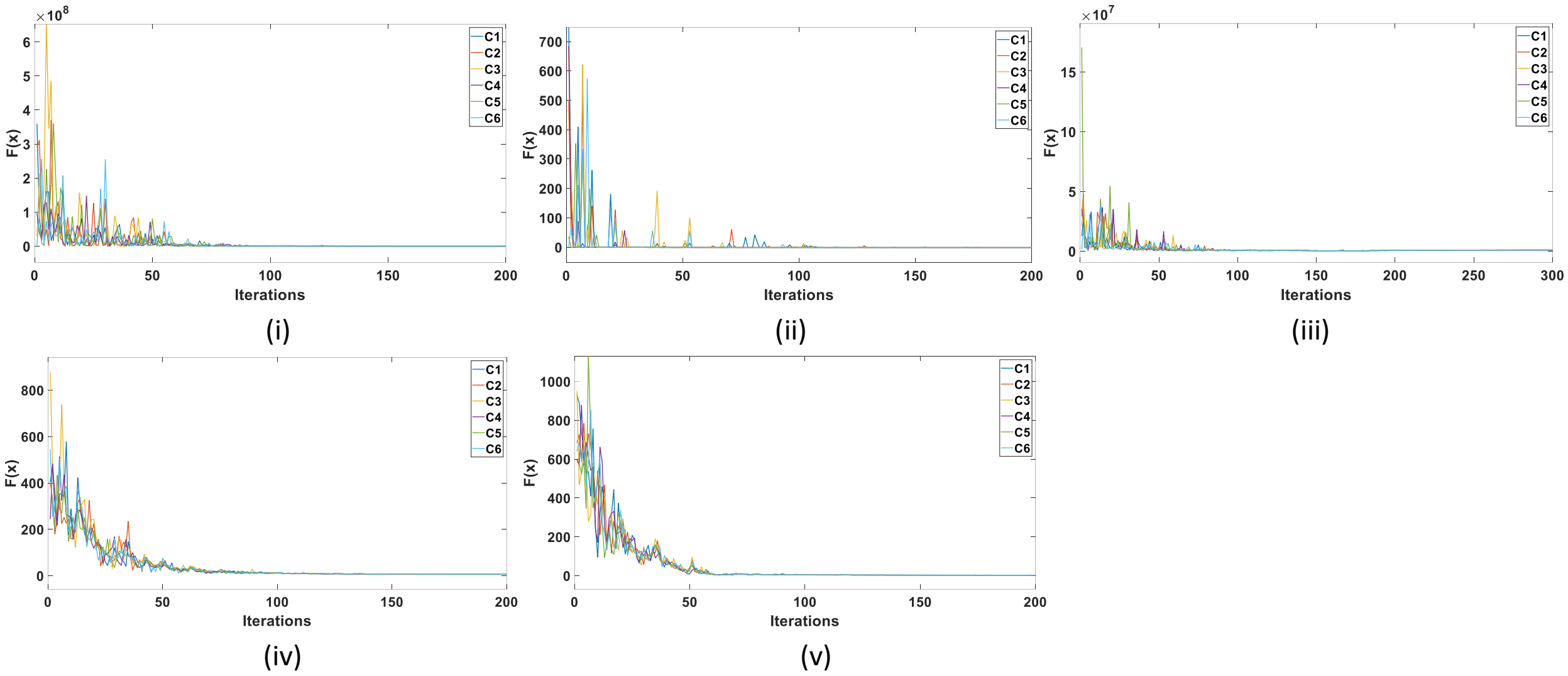

Fig. 4: Convergence plots of Industrial Chemical Process Problems using CI-SAPF-CBO algorithm
(i) Optimal Operation of Alkylation Unit, (ii) Reactor Network Design, (iii) Haverly's Pooling Problems, (iv) Blending Pooling Separation Problem, (v) Propane, Isobutane, n-Butane Nonsharp Separation

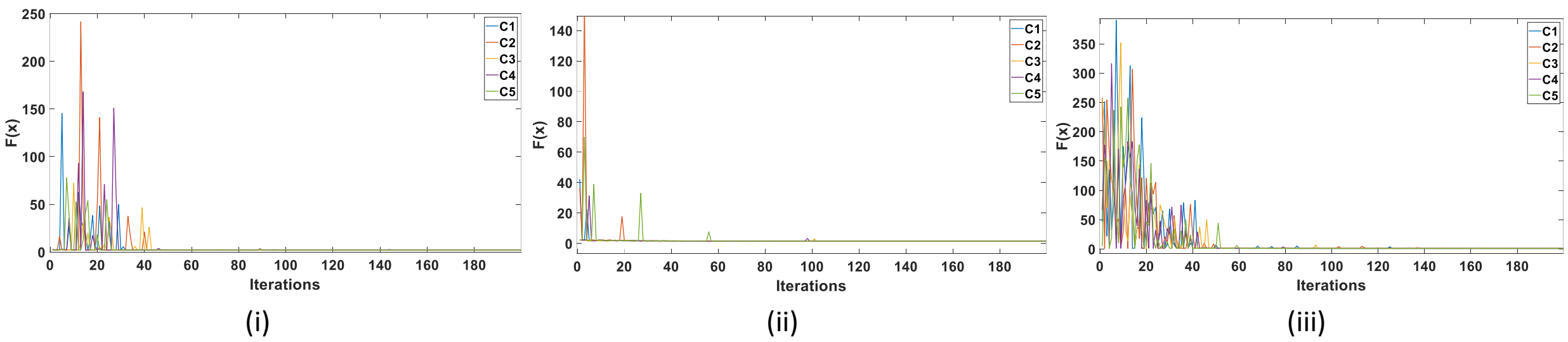

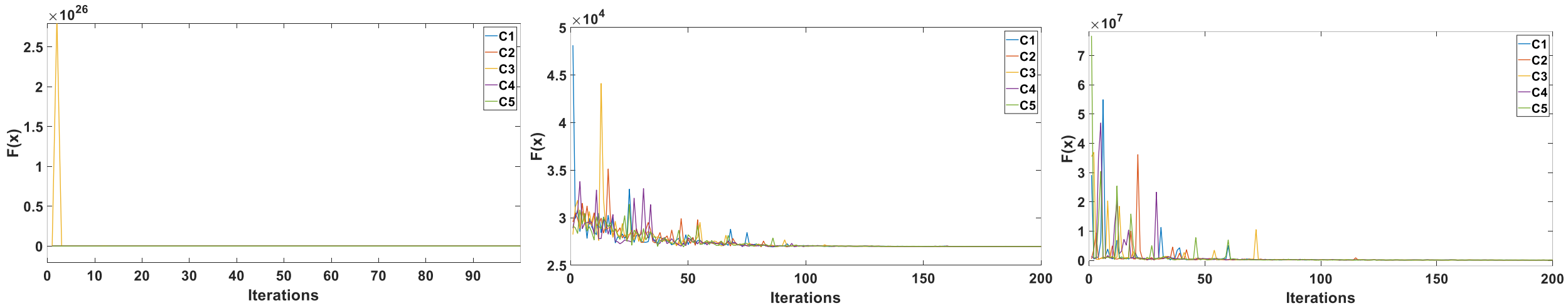

(iv) (v) (vi)

Fig. 5: Convergence plots of Process Synthesis and Design Problems using CI-SAPF algorithm
(i) Process synthesis problem Case I, (ii) Process synthesis and design problem, (iii) Process flow sheeting problem, (iv) Process Synthesis Problem Case 2, (v) Process design Problem, (vi) Multi-product batch plant

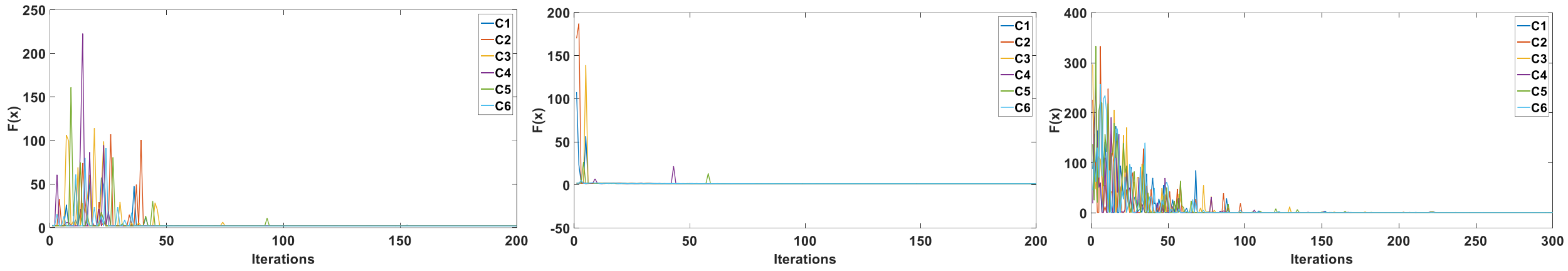

(i) (ii) (iii)

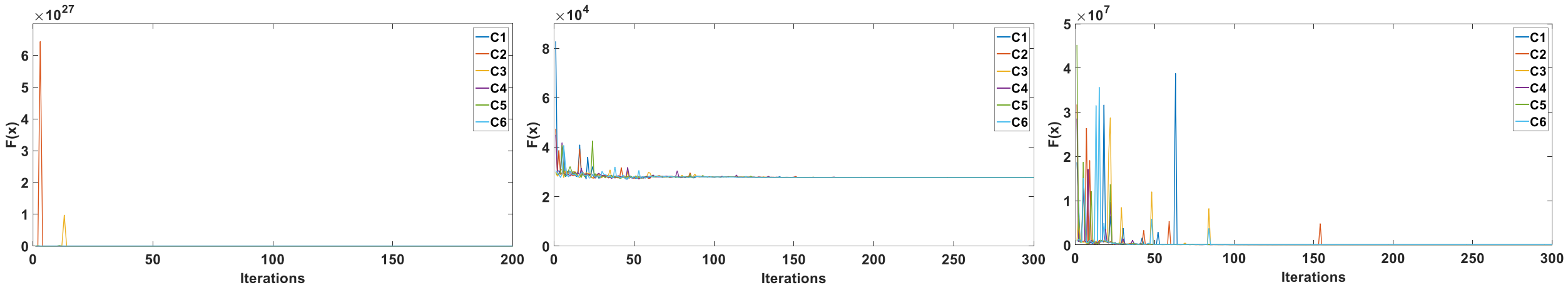

(iv) (v) (vi)

Fig. 6: Convergence plots of Process Synthesis and Design Problems using CI-SAPF-CBO algorithm
(i) Process synthesis problem Case I, (ii) Process synthesis and design problem, (iii) Process flow sheeting problem, (iv) Process Synthesis Problem Case 2, (v) Process design Problem, (vi) Multi-product batch plant

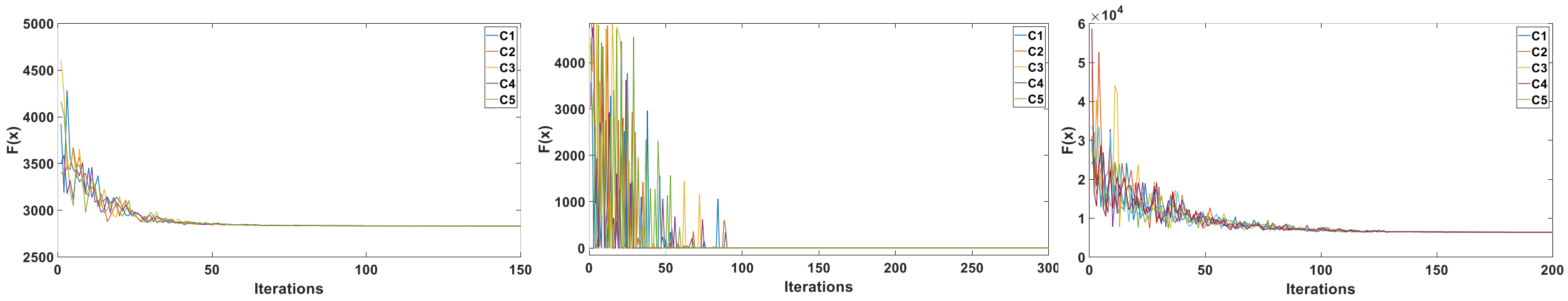

(i) (ii) (iii)

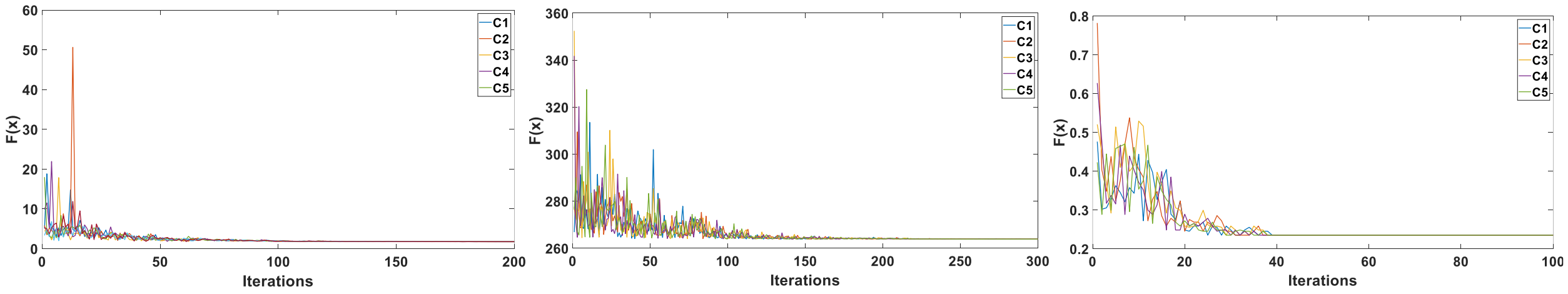

(iv) (v) (vi)

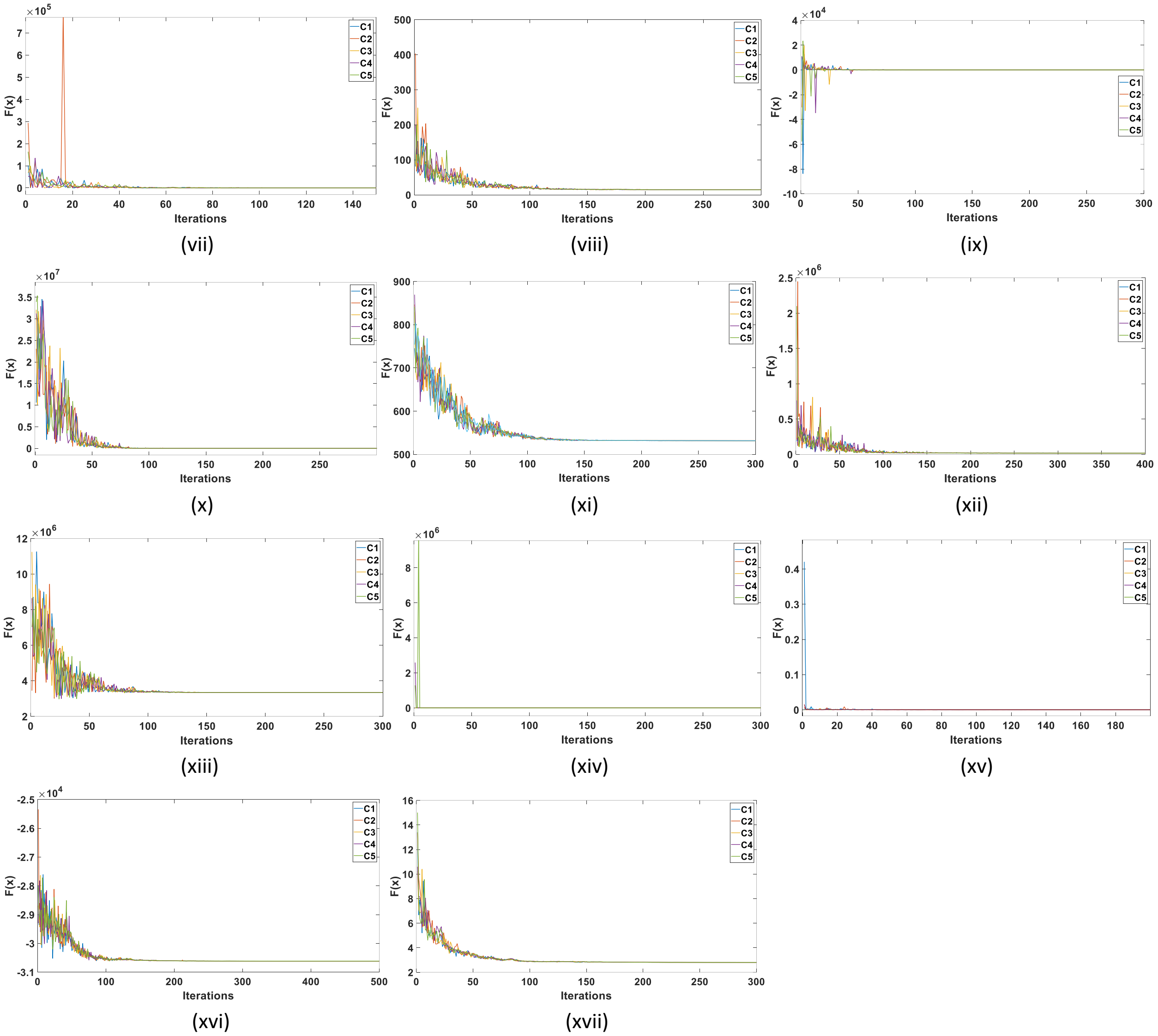

Fig. 8: Convergence plots of Mechanical Engineering Problems using CI-SAPF algorithm
(i) Weight minimization of a speed reducer problem, (ii) Tension/compression spring design Case I, (iii) Pressure vessel design problem, (iv) Welded beam design problem (v) Three bar truss design problem, (vi) Multi disc multi brake design problems, (vii) Planetary gear train design problem, (viii) Step cone pully problem, (ix) Robot gripper problem, (x) Four stage gear box problem, (xi) Ten bar truss design problem, (xii) Rolling element bearing problem, (xiii) Gas transmission compressor design problem, (xiv) Tension/compression spring design Case II, (xv) Gear train design problem, (xvi) Himmelblau's function, (xvii) Topology optimization problem

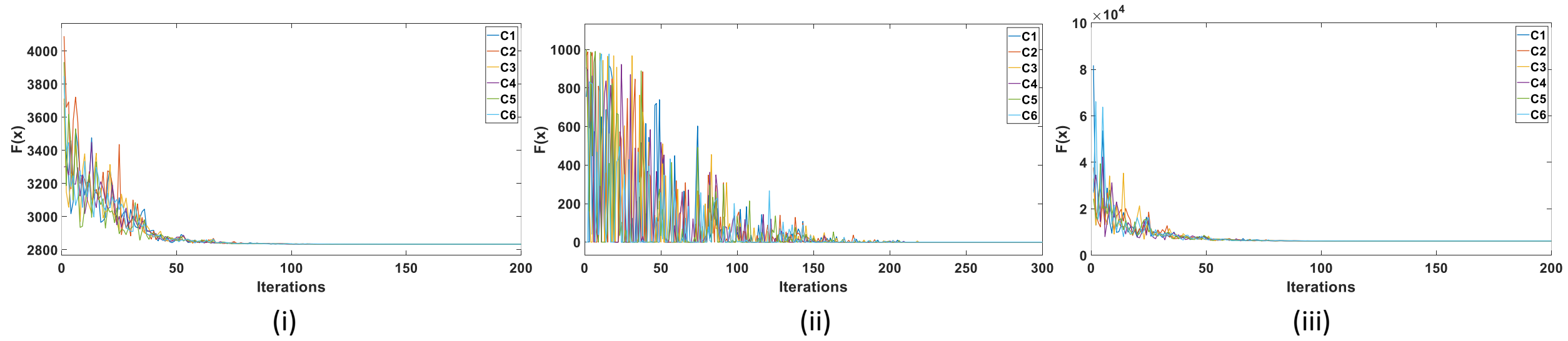

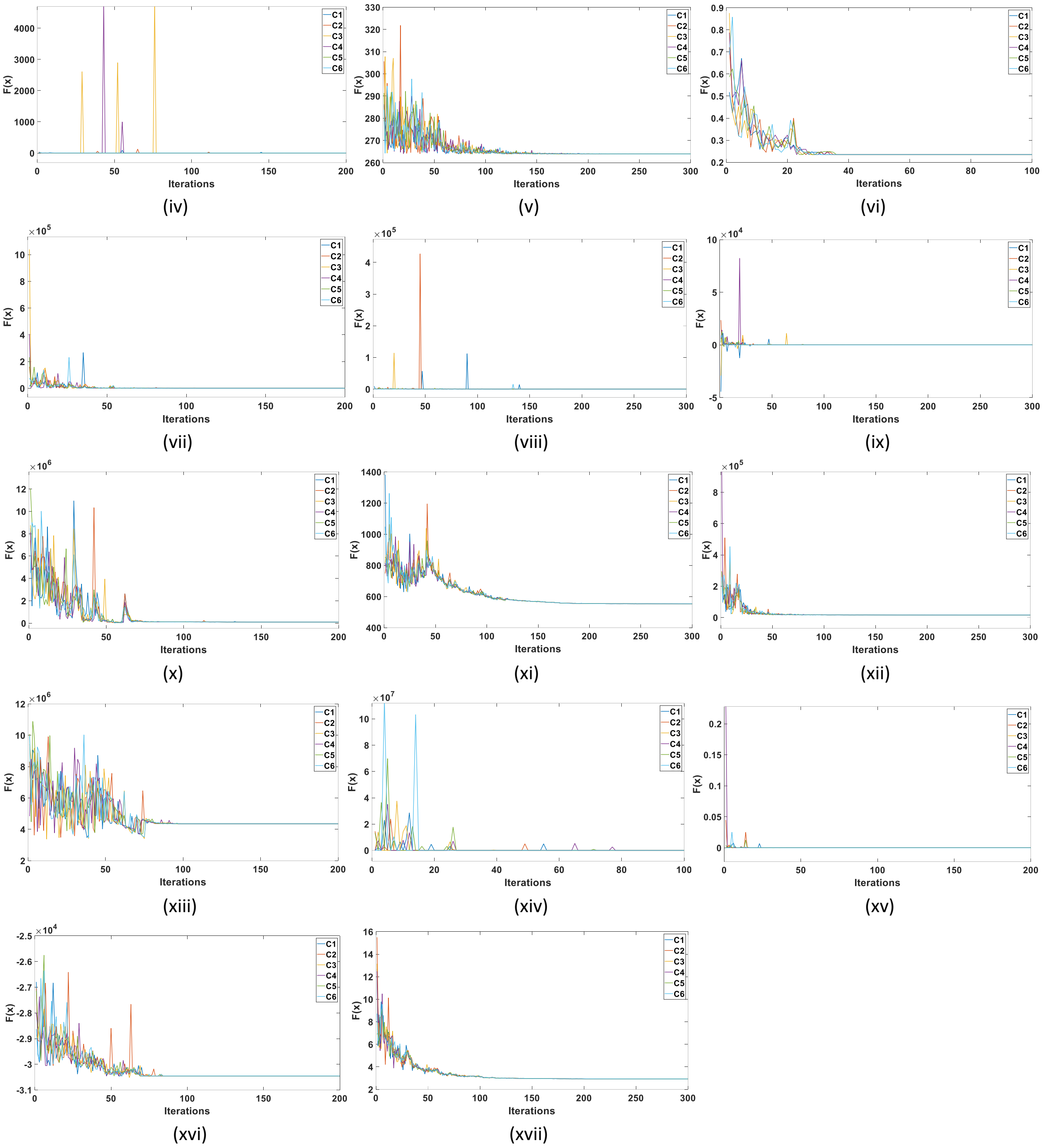

Fig. 7: Convergence plots of Mechanical Engineering Problems using CI-SAPF-CBO algorithm
(i) Weight minimization of a speed reducer problem, (ii) Tension/compression spring design Case I, (iii) Pressure vessel design problem, (iv) Welded beam design problem (v) Three bar truss design problem, (vi) Multi disc multi brake design problems, (vii) Planetary gear train design problem, (viii) Step cone pully problem, (ix) Robot gripper problem, (x) Four stage gear box problem, (xi) Ten bar truss design problem, (xii) Rolling element bearing problem, (xiii) Gas transmission compressor design problem, (xiv) Tension/compression spring design Case II, (xv) Gear train design problem, (xvi) Himmelblau's function, (xvii) Topology optimization problem

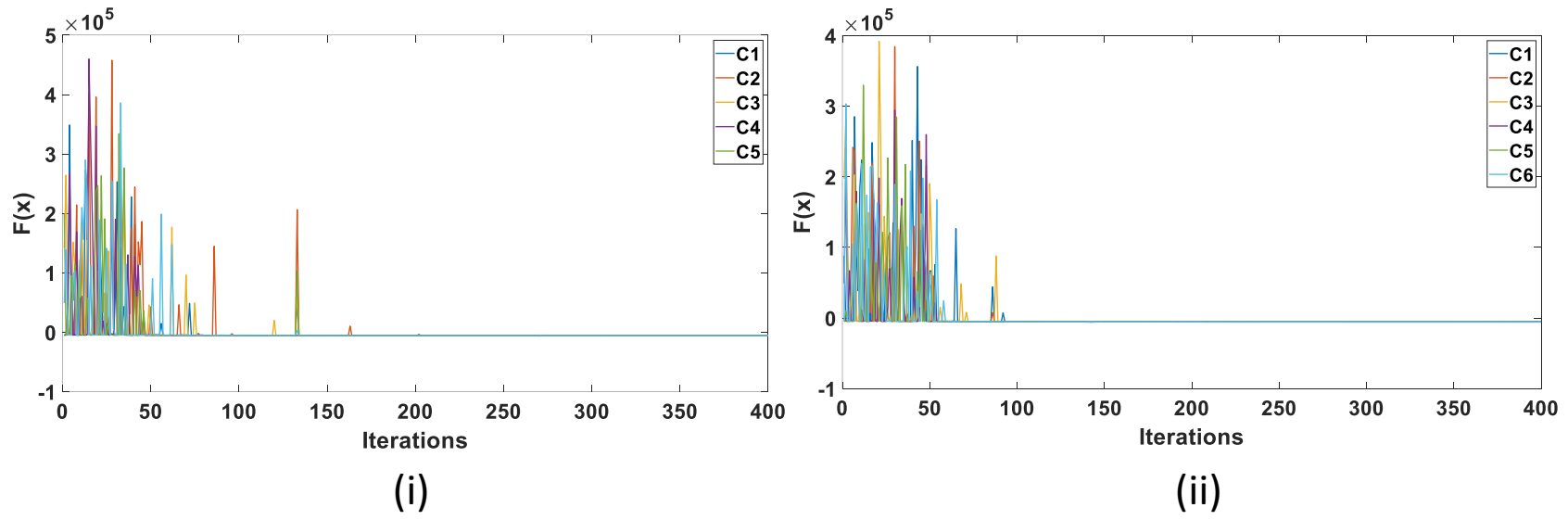

Fig. 9: Convergence plot of wind farm layout problem using (i) CI-SAPF algorithm and (ii) CI-SAPF-CBO algorithm

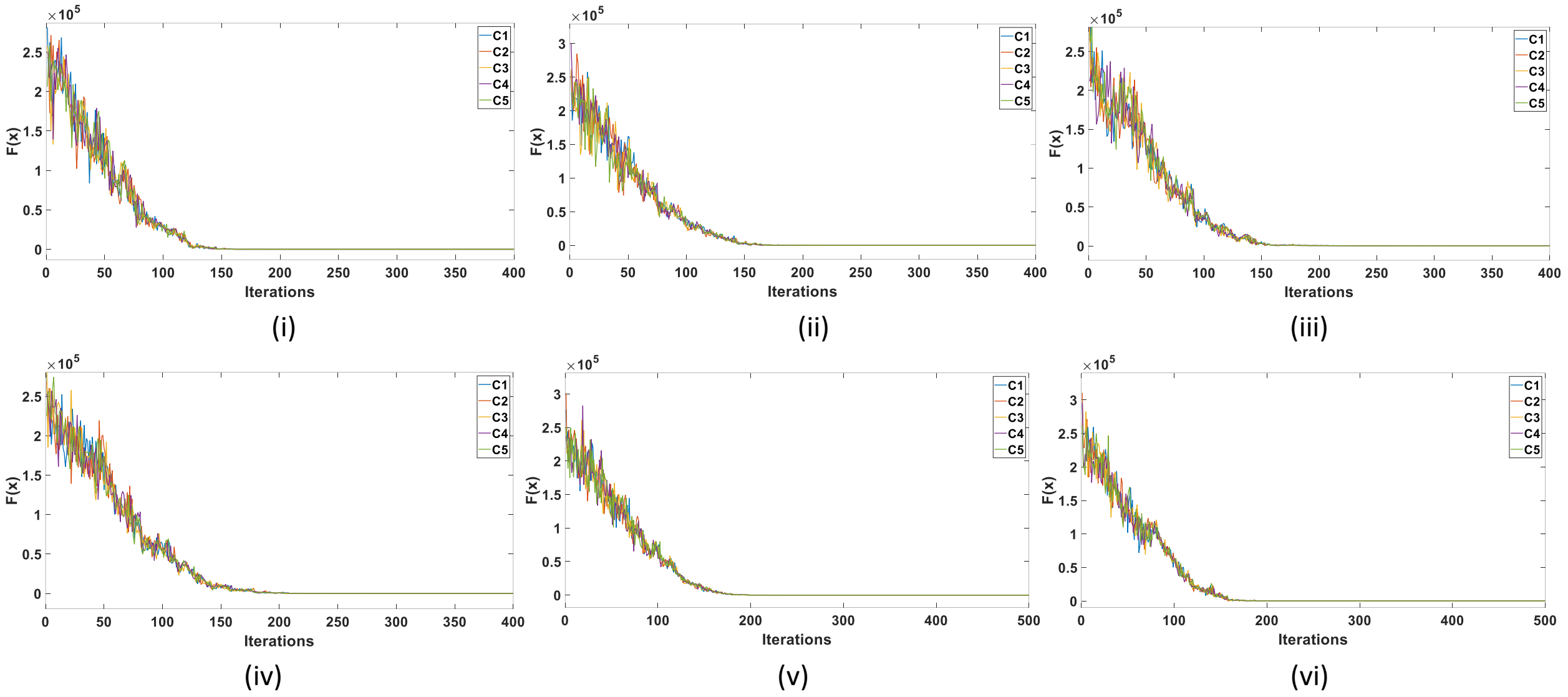

Fig. 10: Convergence plots of power electronic problems using CI-SAPF algorithm
(i) SOPWM for 3-level inverter, (ii) SOPWM for 5-level inverter, (iii) SOPWM for 7-level inverter, (iv) SOPWM for 9-level inverter, (v) SOPWM for 11-level inverter, (vi) SOPWM for 13-level inverter

### 3.2 Result Analysis and Discussion

The constrained optimization test suit problems considered in this work are consisted of both equality and inequality type nonlinear constraints. These problems are from different specialized domains such as industrial chemical process, process synthesis and design, mechanical engineering, power system, power electronic and livestock feed ration (Kumar et al., 2020). These problems were earlier solved using IUDE (Trivedi et al., 2018), $\epsilon$MAg-ES (Hellwig and Beyer, 2018) and iLSHADE$\epsilon$ (Fan et al., 2018) algorithms by Kumar et al. (2020). In this work, CI-SAPF and CI-SAPF-CBO algorithms (Kale and Kulkarni, 2021) are applied to solve 57 constrained optimization problems; out of those, 35 problems are successfully solved using CI-SAPF and 29 problems using CI-SAPF-CBO; however, for rest of the problems they are unable to get the feasible solution. In Kumar et al. (2021), IUDE, $\epsilon$MAg-ES and iLSHADE$\epsilon$ were applied to solve these problems they received the feasible solution for 30, 29 and 35 problems, respectively (refer Table 3). For every problem, the CI-SAPF and CI-SAPF-CBO algorithms were run for 30 times. The result comparison with the best function value available in the literature is presented in Table 1. The statistical results such as best, median, mean and worst function value, standard deviation, function evaluations and computational time obtained using CI-SAPF and CI-SAPF-CBO are compared with IUDE, $\epsilon$MAg-ES and iLSHADE$\epsilon$ are presented in Table 2. For IUDE, $\epsilon$MAg-ES and iLSHADE$\epsilon$ algorithms the function evaluations and computational time are not available in the literature. In Table 2, for CI-SAPF and CI-SAPF-CBO only feasible solutions are shown. To compare the infeasibility, the best, mean and worst constraint violations for infeasible solution are presented and compared with IUDE, $\epsilon$MAg-ES and iLSHADE$\epsilon$ algorithms. The CI-SAPF and iLSHADE$\epsilon$ algorithms obtained the feasible solution for same number of problems i.e., 35. The convergence plots of the solved problems using CI-SAPF and CI-SAPF-CBO algorithms are presented in Figure 3 to Figure 9. Figures 3, 5, 7, 9(i) and 10 represents the convergence plots obtained from CI-SAPF algorithm and Figures 4, 6, 8 and 9(ii) represents the convergence plots obtained from CI-SAPF-CBO for solving industrial chemical process problems, process synthesis and design problems, mechanical engineering problems, power system problem and power electronic problems respectively. In the convergence plots, C1-C5 represents the number of cohort candidates considered the computation. The number candidates considered for CI-SAPF are five and number of candidates for CI-SAPF-CBO are six. It is necessary to consider even number of candidates for CI-SAPF-CBO algorithm as in its mechanism the total number of candidates are divided in two, i.e., slow learning and fast learning candidates.

Table 4: Performance of CI-SAPF and CI-SAPF-CBO as compared to IUDE, $\epsilon$MAgES and iLSHADE$\epsilon$

| **Algorithms** | **CI-SAPF** | | | **CI-SAPF-CBO** | | |
|---|---|---|---|---|---|---|
| | **Winner for best function value** | **Winner for Mean function value** | **Worst function value** | **Winner for best function value** | **Winner for Mean function value** | **Worst function value** |
| **IUDE** | 19 | 17 | 4 | 18 | 15 | 9 |
| **$\epsilon$MAg-ES** | 13 | 15 | 9 | 12 | 12 | 13 |
| **iLSHADE$\epsilon$** | 17 | 16 | 9 | 16 | 12 | 8 |

From the comparison (refer Table 2 and Table 3) it is noticed that CI-SAPF and CI-SAPF-CBO have shown precisely similar performance as compared to other contemporary algorithms discussed in the literature. The winning performance of CI-SAPF and CI-SAPF-CBO as compared to IUDE, $\epsilon$MAg-ES and iLSHADE$\epsilon$ is presented in Table 4. For obtaining the best function value CI-SAPF won 19, 13 and 17 times, and CI-SAPF-CBO won 18, 12 and 16 times as compared to IUDE, $\epsilon$MAg-ES and iLSHADE$\epsilon$, respectively. For mean function value CI-SAPF won 17, 15 and 16 times, and CI-SAPF-CBO won 15, 12 and 12 times as compared to IUDE, $\epsilon$MAg-ES and iLSHADE$\epsilon$, respectively (refer Table 4). However, CI-SAPF and CI-SAPF-CBO considered to be inefficient to handle the power system problems, live-stock feed ration optimization problems and some specific problems having equality constraints. These techniques are superior to solve the problems having inequality constraints.

It is important to discuss the mechanism of the comparative algorithms to execute the rationales behind those mechanisms. In IUDE (Trivedi et al., 2018), the dual population-based approach was modeled incorporated with three mutation strategies such as current-to-pbest, current-to-rand and rand mutation. In the first half of the population, three trail vectors combined with a control parameter to generate a new trial vector. For second half of the population, trail vectors are periodically self-adapted from their own experience which helps to evolve the first half of the population. An improved mutation operator was incorporated in iLSHADE$\epsilon$ (Fan et al., 2018) algorithm which results to evolve the solution quality. In CI-SAPF algorithm, the basic mechanism of CI algorithm is used incorporated with constraint handling SAPF approach. The quality of CI algorithm driven by the self-supervised learning candidates having inherently common goal to achieve the best possible behavior by following and competing with other candidates in the cohort using roulette wheel approach. This mechanism of CI assists to evolve the solution quality to the extent. In CI algorithm, an important parameter i.e., sampling space reduction factor $R$ is in fact drive the solution quality of CI, and it is also one of the limitations which needs several preliminary trials to set an appropriate value. This limitation is overcome by combining the socio-based and physics-based features CI and CBO algorithm, respectively which forms a hybrid CI-SAPF-CBO. In CI-SAPF-CBO algorithm, CI assists to obtain the global optimum and CBO updates the solution iteratively. These are the key feature where CI-SAPF and CI-SAPF-CBO considered to be better than other contemporary algorithms.

As presented in Table 3, $\epsilon$MAg-ES (Hellwig and Beyer, 2018) algorithm obtained the feasible solution for 29 number of problems this is due to the gradient-based approach requires additional function evaluations per execution. This seems to be computationally inefficient and does not guarantee the feasible solution. Further, $\epsilon$MAg-ES requires various computational parameters that need to be readjusted/ tuned which may deteriorate the solution quality and need to invest several preliminary trials to set an appropriate value. The CI-SAPF-CBO also found the feasible solution for 29 problems however, it has obtained better results for 12 problems as compared to $\epsilon$MAg-ES (refer Table 4). In the mechanism of CI-SAPF-CBO, the fast-learning candidate promotes the slow learning candidate to achieve the better position in the search space. Due to which results to explore search space and candidates get better chance to follow the other candidates.

The constrained handling approaches, in IUDE a combined approach of feasibility-based rule and $\epsilon$-constraint were used. In iLSHADE$\epsilon$, an improved $\epsilon$ (I$\epsilon$) constraint handling approach in which improved $\epsilon$ was combined with feasibility-based approach. In I$\epsilon$ constraint handling approach, the value of $\epsilon$ is adaptively adjusted based on the proportion of current feasible solutions which results to balance the search between feasible and infeasible region. In $\epsilon$MAg-ES two constraint handling approaches has been introduced i.e., $\epsilon$ level ordering and a gradient-based repair step approach. However, as discussed earlier it does not guarantee the feasible solution. In CI-SAPF and CI-SAPF-CBO algorithms, the SAPF approach plays a key role in constraint handling. SAPF is a new adaptive version of dynamic penalty function approach. The iteratively generated function value (behavior) $f(\boldsymbol{X})$ for every individual candidate itself work as a penalty parameter which shows the decentralized approach towards constraint handling. This approach is free from parameter tuning so it does not require any preliminary trials to set an appropriate value/s.

## 5. Conclusion and future directions

Nature inspired optimization techniques have been developed to ease the problems solving process. Initially, newly developed algorithms are tested on various unconstrained and constrained benchmark functions. However, it is very much important to validate the efficacy of those algorithm to solve complex real-world discrete and mixed variable constrained optimization problems. In this work, the CI-SAPF and CI-SAPF-CBO algorithms are applied to solve more complex real-world constrained optimization problems. These problems are associated with linear, nonlinear, convex and non-convex functions and similar type of equality and inequality constraints. Among a set of 57 problems, 35 and 29 problems are successfully solved using CI-SAPF and CI-SAPF-CBO, respectively. For rest of the problems infeasible solutions are obtained. Whereas, the IUDE, $\epsilon$MAg-ES and iLSHADE$\epsilon$ algorithms received feasible solutions for 30, 29 and 35 problems respectively. From the solved problems CI-SAPF and CI-SAPF-CBO has shown better performance in all the terms such as feasibility, function value and standard deviations as compared to IUDE, $\epsilon$MAg-ES and iLSHADE$\epsilon$. In CI-SAPF, an intrinsic property of CI i.e., roulette wheel approach supervised the cohort candidates to learn from the other candidate and achieve the best possible behavior and sampling parameter $(R)$ drives the solution quality. In CI-SAPF-CBO, the mechanism of CBO associated with the $COR$ helps to refined the solution obtained from CI-SAPF without parament $R$. It has been noticed that the constraint handling SAPF approach can efficiently handle large number of inequality constraints; however, it does not consider to be more efficient where equality constraints are available in large numbers. In particular, CI-SAPF and CI-SAPF-CBO have shown a descent performance towards handling more complex real-world applications. Though, it is tedious to deal with equality constraints, it necessitates to propose an efficient constraint handling technique in order to handle large number of equality and inequality constraints.